\documentclass[11pt]{article}
\usepackage[utf8]{inputenc}
\usepackage{amsmath}
\usepackage{amsfonts}
\usepackage{graphicx}
\usepackage{subcaption}
\usepackage{fullpage}
\usepackage{todonotes}
\usepackage{physics}
\usepackage{bm}
\usepackage{amssymb}

\newtheorem{remark}{Remark}

\graphicspath{ {/Users/zioepa/PycharmProjects/Cemracs /project_cemracs/Data/} }

\title{CEMRACS 2017: Numerical Probabilistic Approach to MFG}
\author{Andrea Angiuli\thanks{Department of Statistics and Applied Probability, University of California Santa Barbara} \and Christy V. Graves\thanks{Program in Applied and Computational Mathematics, Princeton University} \and  Houzhi Li\thanks{Laboratoire de Probabilités, Modélisation et Statistiques, Université Paris Diderot} \and Jean-Fran\c{c}ois Chassagneux$^\ddagger$ \and Fran\c{c}ois Delarue\thanks{Laboratoire Jean-Alexandre Dieudonné, Université de Nice Sophia-Antipolis} \and René Carmona\thanks{Operations Research and Financial Engineering, Bendheim Center for Finance, Princeton University}}
\date{\today}

\begin{document}

\maketitle

\begin{abstract}
This project investigates numerical methods for solving fully coupled forward-backward stochastic differential equations (FBSDEs) of McKean-Vlasov type. Having numerical solvers for such mean field FBSDEs is of interest because of the potential application of these equations to 
optimization problems over a large population, say for instance mean field games (MFG) and optimal mean field control problems. 
Theory for this kind of problems has met with great success since the early works on mean field games by Lasry and Lions, see \cite{Lasry_Lions}, and by Huang, Caines, and Malham\'{e}, see \cite{Huang}. Generally speaking, the purpose 
is to understand the continuum limit of optimizers or of equilibria (say in Nash sense) as the number of underlying players tends to infinity. When approached from the probabilistic viewpoint, solutions to these control problems (or games) can be described by coupled mean field FBSDEs, meaning that the coefficients depend upon the own marginal laws of the solution. In this note, we detail two methods for solving such FBSDEs which we implement and apply to five benchmark problems. The first method uses a tree structure to represent the pathwise laws of the solution, whereas the second method uses a grid discretization to represent the time marginal laws of the solutions. Both are based on a Picard scheme; importantly, we combine each of them with a generic continuation method that permits to extend the time horizon (or equivalently the coupling strength between the two equations) for which the Picard iteration converges.

\end{abstract}

\section{Introduction}
\hspace{5mm} In this project, we examine numerical methods for solving forward backward stochastic differential equations (FBSDEs) of McKean-Vlasov type. We are particularly interested in equations of this type as they can be used to represent, from the probabilistic viewpoint, solutions to mean field games or, more generally, to mean field stochastic control problems. 

Mean field games were developed independently and at about the same time by Lasry and Lions \cite{Lasry_Lions}, and Huang, Caines, and Malham\'{e} \cite{Huang}. The goal of this theory is to understand the limit as $N \rightarrow \infty$ of the Nash equilibria for an $N$ player stochastic dynamic game with mean field interaction, meaning that players interact with one another through their collective state.  Equivalently, mean field games must be regarded as the continuum limit of games with a large number of symmetric players, each 
of them having a small effect on the dynamics of the whole group. The applications of mean field games are numerous, and spread across many disciplines, including social science (congestion \cite{congestion1}\cite{congestion2}\cite{congestion3}, cyber attacks \cite{cyber}), biology (flocking \cite{Nourian}\cite{flocking2}), and economics (systemic risk \cite{systemic_risk}, production of exhaustible resources \cite{sircar} \cite{chan}), just to name a few.
As explained in
 \cite{MFGP,FBSDEs}, solutions to mean field games can be characterized through a coupled system of two forward and backward stochastic differential equations of mean field type, like those we address in this note. The forward equation provides the dynamics of one typical player in the population at equilibrium. Generally speaking, the backward equation accounts for the optimality condition in the definition of an equilibrium and the McKean-Vlasov structure is precisely here to stress the fact that the player in hand is representative of the others.  
As exemplified in 
\cite{CDL}, other forms of equilibria can be addressed by means of this kind of equations. This includes optimal mean field control problems, which can be regarded as the continuum limit of a control problem involving a large number of symmetric players obeying a central planner. Below, we mostly focus on examples arising in the theory of mean field games.  

Whilst deterministic numerical methods, based upon finite differences or variational approaches, are also conceivable for handling mean field games, see 
\cite{AchdouCapuzzo-Dolcetta,AchdouCamilliCapuzzo-Dolcetta2,AchdouPorretta}
and 
\cite{BenamouCarlier,LachapelleSalomonTurinici,Gueant_numerical}, we here focus on the approach based on FBSDEs. In this regard, we
implement (and compare)  two different algorithms. The first algorithm, which is based on the paper of Chassagneux, Crisan, and Delarue \cite{Solver}, relies on a tree structure to represent the pathwise law of the solution. The second algorithm, and main contribution of this report, takes the algorithm presented in the paper of Delarue and Menozzi \cite{Grid} for solving FBSDEs and extends it to the mean field framework in hand. In this algorithm, a grid structure is used to represent the marginal laws of the solution. 
The serious issue that we are facing in this note is that both methods are based upon a Picard scheme, the first method involving a global Picard scheme upon the whole process and the second one involving a Picard scheme on the sole marginal laws of the process. 
It is indeed a well-known fact that, because of the strong coupling between the forward and backward equations, Picard schemes for FBSDEs may just converge in small time, even in the classical case without mean field interaction. For sure, this limitation should persist 
in the mean field setting for the global Picard method; as exemplified below, it turns out that it persists as well when the Picard scheme is applied to the marginal laws.
One of our main contribution in this report is to apply the 
time continuation approach presented in \cite{Solver} to the grid algorithm and to compare the results with the tree algorithm for which the time continuation approach was originally designed in \cite{Solver}. In brief, the time continuation permits to extend, by a continuation argument, the time interval on which the Picard scheme converges. We illustrate both algorithms on a handful of example problems.

Section \ref{sec:review} provides a review of Nash equilibria in $N$ player stochastic differential games, and their continuum mean field game counterparts. We review two probabilistic approaches to formulate the solutions of mean field games and provide the general FBSDE system which we would like to solve. In Section \ref{Algorithms}, we describe the algorithms that we implement in the report. Some benchmark examples and the corresponding numerical results are presented in Section \ref{Examples}. We conclude in Section \ref{Conclusion}.

\section{Overview of Mean Field Games and FBSDEs}\label{sec:review}

The purpose of this section is to introduce the theoretical material that is needed for our numerical analysis. 
The objective is purely pedagogical and the text does not contain any new result.

\subsection{$N$ Player Stochastic Differential Games} \label{N-Player}
\hspace{5mm} 
We start with the description of the prototype of a finite player game in the theory of mean field games. We consider $N \in \mathbb{Z}^+$ players indexed by $i \in \{1,\dots, N \}$. The dynamic game occurs over a fixed time horizon $[0,T]$ for some $T>0$. We have $N$ independent $m$-dimensional Brownian motions $(W^i_t)_{0\leq t \leq T}$ which are supported by a filtered probability space $(\Omega,\mathcal{F},\mathbb{F}=(\mathcal{F}_t)_{0 \leq t \leq T},\mathbb{P})$. Each player chooses its control $\alpha^i=(\alpha^i_t)_{0\leq t\leq T}$ from the set $\mathbb{A}$ defined as the set of square integrable $\mathbb{F}$ adapted processes with values in a given set $A$ (typically $A$ is a closed convex subset of a Euclidean space). Each player $i$ has a state $X^i$ which evolves according to the stochastic differential equation:

\begin{equation*}
    dX^i_t=b^i(t,X^i_t,\bar{\mu}_t,\alpha_t)dt+\sigma^i(t,X^i_t,\bar{\mu}_t,\alpha_t) dW^i_t,
\end{equation*}
where $\bar{\mu}_t$ denotes the empirical distribution of the players' states: $\bar{\mu}_t=\frac{1}{N}\sum_{j=0}^N \delta_{X^j_t} \in \mathcal{P}_2(\mathbb{R}^d)$. Here $\mathcal{P}_2(\mathbb{R}^d)$ is the space of probability measures with a finite second moment, which we equip with the 2-Wasserstein distance, denoted by $W_2$. For $\mu, \nu \in \mathcal{P}_2(\mathbb{R}^d)$, we call $\Gamma(\mu, \nu)$ the
set of all the joint laws with marginals $\mu$ and $\nu$. Then, the 2-Wasserstein distance is defined by:

\begin{equation*}
 W_2(\mu,\nu) = \left(\inf_{\gamma \in \Gamma(\mu, \nu)} \int |x - y|^2 d \gamma(x,y)\right)^{1/2}.
\end{equation*}

The drift and volatility functions, $b^i$ and $\sigma^i$, respectively, are deterministic functions $(b^i,\sigma^i): [0,T] \cross \mathbb{R}^d \cross \mathcal{P}_2(\mathbb{R}^d) \cross {A} \rightarrow \mathbb{R}^d \cross \mathbb{R}^{d \cross m}$. Most of the time, they are assumed to be bounded in time and to be Lipschitz continuous with respect to all the arguments, the Lipschitz property in the measure argument being understood with respect to $W_{2}$. This ensures that, for a given $\bm{\alpha}=(\alpha^1,\cdots,\alpha^N)$, the state dynamics 
$(X^1,\cdots,X^N)$ is well defined. 

Given a tuple of controls $\bm{\alpha}=(\alpha^1,\dots \alpha^N)$, we associate with player $i$ a cost objective which we take to be of the form:

\begin{equation*}
    J^i(\bm{\alpha})=\mathbb{E} \left[\int_0^T f^i(t,X^i_t,\bar{\mu}_t,\alpha^i_t) dt+ g^i(X^i_T,\bar{\mu}_T) \right].
\end{equation*}

Thus each player considers a deterministic running cost $f^i:[0,T] \cross \mathbb{R}^d \cross \mathcal{P}_2(\mathbb{R}^d) \cross {A} \rightarrow \mathbb{R}$, and deterministic terminal cost $g^i:\mathbb{R}^d \cross \mathcal{P}_2(\mathbb{R}^d) \rightarrow \mathbb{R}$. Of course, each of them wishes to minimize its own cost by tuning its own control in the most relevant way. 
Note that we only allow the interaction of the players through their empirical measure, as this will be needed in our formulation of the continuum limit. Still, extensions exist, in which players also interact through the controls, see Subsection \ref{extended}.

The players are in a \textit{Nash equilibrium} if each player is no better off for switching their strategy when they consider the other players' strategies to be fixed. More precisely, the set of strategies $\bm{\alpha}$ is a Nash equilibrium if

\begin{equation*}
    J^i(\bm{\alpha}) \leq J^i(\alpha^1,\dots, \alpha^{i-1},\alpha,\alpha^{i+1},\dots, \alpha^N), \forall \alpha \in \mathbb{A}, \forall i \in \{1,\dots,N\}.
\end{equation*}

\subsection{Mean Field Games} \label{MFG}
\hspace{5mm} For games where $N$ is large, the problem quickly becomes of an intractable complexity. Thus we turn to the continuum limit by considering the limit as $N$ tends to infinity. 
In order for this limit to make sense, we require the players to be symmetric. Precisely, we require $b=b^i$, $\sigma=\sigma^i$, $f=f^i$, and $g=g^i$ $\forall i \in \{1, \dots, N \}$. As the number of players increase, the impact of each player on the empirical distribution decreases, and we expect to have a propagation of chaos such that the players become asymptotically independent of each other. This is the rationale for passing to the limit: Asymptotically, 
the influence of one player on the group should be null
and 
the statistical structure of the whole should be pretty simple.

We wish to formulate the analogue of a Nash equilibrium when there is a continuum of players. To this end, we consider the states and actions of the other players to be fixed, and consider the best response for a representative player (as we expect equilibria to inherit the symmetric structure of the game). Thus, the first step is to solve an optimization problem. The next step is to find a fixed point, providing an analogue of a Nash equilibrium for the mean field game.

We again have a filtered probability space $(\Omega,\mathcal{F},\mathbb{F}=(\mathcal{F}_t)_{0 \leq t \leq T},\mathbb{P})$ where the filtration supports an $m$-dimensional Brownian motion $W=(W_t)_{0\leq t \leq T}$ and an initial condition $\xi \in L^2(\Omega,\mathcal{F}_0,\mathbb{P};\mathbb{R}^d)$. 
 
The strategy for solving the asymptotic game is the following:

\begin{enumerate}
    \item For a fixed deterministic flow of probability measures $\mu=(\mu_t)_{0\leq t \leq T} \in \mathcal{C}([0,T] , \mathcal{P}_2(\mathbb{R}^d))$, solve the standard stochastic control problem:
    \begin{equation}
    \label{SSCP}
    \inf_{\alpha \in \mathbb{A}}J^{\mu}(\alpha)=\mathbb{E}\left[\int_0^T f(t,X^{\alpha}_t,\mu_t,\alpha_t)dt+g(X^{\alpha}_t,\mu_T) \right],
    \end{equation}
    subject to
    
    \begin{equation*}
    \begin{split}
    dX^{\alpha}_t&=b(t,X^{\alpha}_t,\mu_t,\alpha_t)dt+\sigma(t,X^{\alpha}_t,\mu_t,\alpha_t) dW_t \\
    X^{\alpha}_0&=\xi.
    \end{split}
    \end{equation*}
    
    \item Find a fixed point, $\mu$, such that $\mathcal{L}(X^{\alpha}_t)=\mu_t$ for all $0\leq t \leq T$.
\end{enumerate}

This strategy can be tackled from either the PDE viewpoint (leading to a coupled Hamilton-Jacobi-Bellman and Kolmogorov/Fokker-Plank equations, known as the MFG system in the literature) \cite{Lasry_Lions} \cite{Huang} or the probabilistic viewpoint \cite{MFGP} \cite{FBSDEs}, which is the focus of this project. Within the probabilistic viewpoint, there are two approaches, both of which are formulated with FBSDEs. See Chapters $3$ and $4$ of the manuscript of Carmona and Delarue \cite{MFGP} for reference on the two probabilistic approaches.

For simplicity, from now on we assume $m$, the dimension of the Brownian motion matches $d$, the dimension of the state variable. We also assume the diffusion coefficient, $\sigma$, is a constant matrix $\sigma \in \mathbb{R}^{d \cross d}$. For both approaches, we will utilize the Hamiltonian deriving from the aforementioned stochastic control problem \eqref{SSCP}. In fact, since we assume that the drift is uncontrolled, we can just write the reduced Hamiltonian:

\begin{equation*}
    H(t,x,\mu,\alpha,y)=b(t,x,\mu,\alpha) \cdot y+f(t,x,\mu,\alpha),
\end{equation*}
for $t \in [0,T]$, $x \in \mathbb{R}^d$, $\mu \in \mathcal{P}_2(\mathbb{R}^d)$, $\alpha \in A$, and an adjoint variable $y \in \mathbb{R}^d$. Then, a key object in order to formulate the solution to 
\eqref{SSCP}
is (whenever it exists):

\begin{equation*}
    \hat{\alpha}(t,x,\mu,y)=\arg \inf_{\alpha \in A} H(t,x,\mu,\alpha,y).
\end{equation*}
We will provide below explicit examples for $\hat{\alpha}(t,x,\mu,y)$. We now give a brief introduction to the two probabilistic approaches.

\subsubsection{Weak approach}
\hspace{5mm} In the first probabilistic approach, which we will refer to as the weak approach, the optimization problem is solved using a backward SDE for the probabilistic representation of the value function. For a fixed flow of measures $\mu=(\mu_{t})_{0 \le t \le T}$, let $u:[0,T] \times \mathbb{R} \rightarrow \mathbb{R}$ denote the value function:

\begin{equation*}
    u(t,x):=\inf_{(\alpha_s)_{t \leq s \leq T} \in \mathbb{A}} \mathbb{E} \left[\int_t^T f(s,X_s,\mu_s,\alpha_s) ds+g(X_T,\mu_T) \mid X_t=x \right].
\end{equation*}

The strategy is to evaluate the value function along the solution of the state process $(X_{t})_{0 \leq t \leq T}$, namely we
let $Y_t=u(t,X_t)$. The weak formulation of the stochastic control problem underpinning $u$ says that, under suitable assumptions 
that are exemplified below, 
the pair $(X_{t},Y_{t})_{0 \leq t \leq T}$
has to solve the following FBSDE:
\begin{equation}
\begin{split}
    dX_t&=b\left(t,X_t,\mu_{t},\hat{\alpha}\left(t,X_t,\mu_{t},\sigma^{-1}Z_t\right)\right)dt+\sigma dW_t \\
    X_0&=\xi ,\\
    dY_t&=-f\left(t,X_t,\mu_{t},\hat{\alpha}\left(t,X_t,\mu_{t},\sigma^{-1}Z_t\right)\right)dt+Z_t dW_t \\
    Y_T&=g(X_T,\mu_{T}),
\end{split}
\label{weak0}
\end{equation}
where we assume $\sigma$ to be invertible. For instance, we take the following set of assumptions from Chapter 3 of the manuscript by Carmona and Delarue \cite{MFGP}: We may assume that the set $A$ for the values of the controls is a bounded closed convex subset of a Euclidean space, the deterministic functions $b$, $f$, and $g$ are defined on $[0,T] \cross \mathbb{R}^d \cross \mathcal{P}_2(\mathbb{R}^d) \cross A$, $[0,T] \cross \mathbb{R}^d \cross \mathcal{P}_2(\mathbb{R}^d) \cross A$, and $\mathbb{R}^d \cross \mathcal{P}_2(\mathbb{R}^d)$, respectively, and there exists a constant $C_0>1$ such that:
\begin{itemize}
\item For any $t\in [0,T]$, $x$, $x'$ $\in \mathbb{R}^d$, $\alpha$, $\alpha '$ $\in A$ and $\mu \in \mathcal{P}_2(\mathbb{R}^d)$ :
\begin{equation*}
\begin{split}
    \abs{(b,f)(t,x',\mu,\alpha')-(b,f)(t,x,\mu,\alpha)}&+\abs{\sigma(t,x')-\sigma(t,x)}+\\
    &+\abs{g(x',\mu)-g(x,\mu)}\leq C_0\abs{(x,\alpha)-(x ',\alpha ')}.
\end{split}
\end{equation*}
\item The functions $b$, $f$, $\sigma$ and $g$ are bounded by $C_0$.
\item There exists a function
\begin{equation*}
\hat{\alpha}:[0,T] \cross \mathbb{R}^d \cross \mathcal{P}_2(\mathbb{R}^d) \cross \mathbb{R}^d \ni (t,x,\mu,y) \rightarrow \hat{\alpha}(t,x,\mu,y)
\end{equation*}

which is $C_0$-Lipschitz continuous in $(x,y)$ such that, for each $(t,x,\mu,y) \in [0,T] \cross \mathbb{R}^d \cross \mathcal{P}_2(\mathbb{R}^d) \cross \mathbb{R}^d $, $\hat{\alpha}(t,x,\mu,y)$ is the unique minimizer of the Hamiltonian $H(t,x,\mu,y,\alpha)$.
\end{itemize}

Under this set of assumptions, it is shown in Chapter 3 of the manuscript by Carmona and Delarue \cite{MFGP} that a flow of measures $\mu=(\mu_{t})_{0 \leq t \leq T}$ is a mean field game equilibrium if and only if $\mu_t=\mathcal{L}(X_t), \forall t \in[0,T]$, where $(X,Y,Z)$ is a solution of the weak approach FBSDE system in Equation (\ref{weak0}), in which case 
(\ref{weak0}) becomes
\begin{equation}
\begin{split}
    dX_t&=b\left(t,X_t,\mathcal{L}(X_t),\hat{\alpha}\left(t,X_t,\mathcal{L}(X_t),\sigma^{-1}Z_t\right)\right)dt+\sigma dW_t \\
    X_0&=\xi ,\\
    dY_t&=-f\left(t,X_t,\mathcal{L}(X_t),\hat{\alpha}\left(t,X_t,\mathcal{L}(X_t),\sigma^{-1}Z_t\right)\right)dt+Z_t dW_t \\
    Y_T&=g(X_T,\mathcal{L}(X_T)),
\end{split}
\label{weak}
\end{equation}
where we use the generic notation ${\mathcal L}(\cdot)$ for the law of a random variable.
This approach is developed further in the papers and manuscript of Carmona and Delarue \cite{MFGP}\cite{FBSDEs}.

\subsubsection{Pontryagin approach}
\hspace{5mm} The second probabilistic approach, which we will refer to as the Pontryagin approach, is based on the Pontryagin stochastic maximum principle. In this formulation, the optimization problem is solved using a backward SDE for the probabilistic representation of the spatial derivative of the value function $u$. Formally, 
the strategy is thus to evaluate the process $(X_{t})_{0 \leq t \leq T}$ along 
$\nabla_{x} u$. Hence we let $Y_t=\nabla_x u(t,X_t)$,
which makes sense when $\nabla_{x} u$ is well-defined. In fact the Pontryagin system may be formulated without any further reference to the regularity of $u$, the Pontryagin formulation having the following general form:

\begin{equation}
\begin{split}
    dX_t=&b\left(t,X_t,\mu_t,\hat{\alpha}\left(t,X_t,\mu_t,Y_t \right)\right)dt+\sigma dW_t \\
    X_0=&\xi, \\
    dY_t=&-[\nabla_x b(\left(t,X_t,\mu_t,\hat{\alpha}\left(t,X_t,\mu_t,Y_t \right)\right))\cdot Y_t \\
    &+\nabla_xf\left(t,X_t,\mu_t,\hat{\alpha}\left(t,X_t,\mu_t,Y_t\right)\right)]dt +Z_t dW_t \\
    Y_T=&\nabla_xg(X_T,\mu_{T}),
\end{split}
\label{Pontryagin0}
\end{equation}
where we assume $b$, $f$ and $g$ to be differentiable with respect to $x$. We may use the following set of assumptions taken from Chapter 3 of the manuscript by Carmona and Delarue \cite{MFGP} to guarantee that the Pontryagin system is both a necessary and a sufficient condition of optimality: We assume the coefficients $b$, $f$, and $g$ are defined on $[0,T] \cross \mathbb{R}^d \cross \mathcal{P}_2(\mathbb{R}^d) \cross A$, $[0,T] \cross \mathbb{R}^d \cross \mathcal{P}_2(\mathbb{R}^d) \cross A$, and $\mathbb{R}^d \cross \mathcal{P}_2(\mathbb{R}^d)$, respectively. We also assume that they satisfy:
\begin{itemize}
\item The drift $b$ is an affine function of $(x,\alpha)$ of the form:
\begin{equation*}
b(t, x,\mu,\alpha) =b_0(t,\mu)+b_1(t)x+b_2(t)\alpha,
\end{equation*}
where $b_0: [0,T] \cross \mathcal{P}_2(\mathbb{R}^d) \ni (t,\mu) \mapsto b_0(t,\mu)$, $b_1 : [0,T] \ni t \mapsto b_1(t)$ and $b_2 : [0,T] \ni t \mapsto b_2(t)$ are $\mathbb{R}^d$, $\mathbb{R}^{d \cross d}$ and $\mathbb{R}^{d\cross d}$ valued, respectively, and are measurable and bounded on bounded subsets of their respective domains. 
\item There exist two constants $C_1 >0$ and $C_2 \geq 1$ such that the function $\mathbb{R}^d \cross A \ni (x,\alpha) 
\mapsto f(t,x,\mu,\alpha) \in \mathbb{R}$ is once continuously differentiable with Lipschitz-continuous derivatives (so that $f(t,\cdot,\mu,\cdot)$ is $C^{1,1}$), the Lipschitz constant in $x$ and $\alpha$ being bounded by $C_2$ (so that it is uniform in $t$ and $\mu$). Moreover, it satisfies the following strong form of convexity:
\begin{equation*}
f(t,x ', \mu, \alpha ')-f(t,x, \mu, \alpha)-(x-x',\alpha - \alpha ') \cdot \partial_{(x,\alpha)}f(t,x,\mu,\alpha)\geq C_1 \abs{\alpha ' - \alpha}^2.
\end{equation*}
The notation $\partial_{(x,\alpha)}f$ stands for the gradient in the joint variables $(x,\alpha)$. Finally, $f$, $\partial_x f$ and $\partial_{\alpha}f$ are locally bounded over $[0,T] \cross \mathbb{R}^d \cross \mathcal{P}_2(\mathbb{R}^d) \cross A$.

\item The function $\mathbb{R}^d \cross \mathcal{P}_2(\mathbb{R}^d)  \ni (x,\mu) \mapsto g(x,\mu)$ is locally bounded. Moreover, for any $\mu \in \mathcal{P}_2(\mathbb{R}^d)$, the function $\mathbb{R}^d \ni x \mapsto g(x,\mu)$ is once continuously differentiable and convex, and has a $C_2$-Lipschitz continuous first order derivative.
\end{itemize}

Under these assumptions, it is shown in Chapter 3 of the manuscript by Carmona and Delarue \cite{MFGP} that a flow of measures $\mu=(\mu_{t})_{0 \le t \le T}$ is a mean field game equilibrium if and only if $\mu_t=\mathcal{L}(X_t), \forall t \in[0,T]$, where $(X,Y,Z)$ is a solution of the Pontryagin approach FBSDE system in Equation (\ref{Pontryagin0}), in which case  
(\ref{Pontryagin0}) becomes

\begin{equation}
\begin{split}
    dX_t&=b\left(t,X_t,\mathcal{L}(X_t),\hat{\alpha}\left(t,X_t,\mathcal{L}(X_t),Y_t \right)\right)dt+\sigma dW_t \\
    X_0&=\xi, 
    \\
    dY_t&=-[\nabla_x b(\left(t,X_t,\mathcal{L}(X_t),\hat{\alpha}\left(t,X_t,\mathcal{L}(X_t),Y_t \right)\right))\cdot Y_t \\
    &+\nabla_xf\left(t,X_t,\mathcal{L}(X_t),\hat{\alpha}\left(t,X_t,\mathcal{L}(X_t),Y_t\right)\right)]dt +Z_t dW_t \\
    Y_T&=\nabla_xg(X_T,\mathcal{L}(X_T)).
\end{split}
\label{Pontryagin}
\end{equation}
This approach is also developed further in the papers and manuscript of Carmona and Delarue \cite{MFGP}\cite{FBSDEs}.

\subsubsection{Mean Field Games of Control}
\label{extended}
\hspace{5mm} In many applications, individuals may interact through their controls, instead of their states. One example is an application of trade crowding which was tackled with a mean field game approach in the paper of Cardaliaguet and Lehalle \cite{trader}. There is also a model for price impact in the book of Carmona and Delarue which we take as one of our example problems in Section \ref{trader}. Mean field games where players interact through the law of their controls is sometimes referred to as extended mean field games \cite{extended}, see also Chapter 4 in \cite{MFGP}. We are interested in testing our numerical methods on a certain class of mean field games of control: those in which the interaction is through the marginal distributions, $\mathcal{L}(X_t)$ and $\mathcal{L}(\alpha_t)$. To design our algorithms to handle such a general framework of mean field interaction, we study numerical methods for solving a general FBSDE system which includes the two approaches detailed above, as well as this class of mean field games of control.

\subsection{General System}
\hspace{5mm} We can address both probabilistic formulations for mean field games and a class of mean field games of control simultaneously by considering the following general FBSDE system. We now take the dimension of the state space to be $d=1$ since, for simplicity, our algorithms will be just applied to this case. Let $[X]=\mathcal{L}(X)$ denote the law of a process $X$. With an abuse of notation, let $[X,Y,Z]:=(\mathcal{L}(X),\mathcal{L}(Y),\mathcal{L}(Z))$ denote the laws of the individual processes (and not their joint law). The general system is the following:
\begin{equation}
    \begin{split}
        dX_t &= B(t,X_t,Y_t,Z_t,[X_t,Y_t,Z_t])dt + \sigma dW_t \\
        X_0&=\xi \in L^2(\Omega,\mathcal{F}_0,\mathbb{P};\mathbb{R}), \\
        dY_t &= -F(t,X_t,Y_t,Z_t,[X_t,Y_t,Z_t])dt + Z_t dW_t \\
        Y_T &= G(X_T,[X_T]). \\
    \end{split}
    \label{general_equation}
\end{equation}
The assumption that the coefficients depend, at most, upon the marginal laws ${\mathcal L}(X_{t})$, ${\mathcal L}(Y_{t})$ and ${\mathcal L}(Z_{t})$ of the triple $(X_{t},Y_{t},Z_{t})$ (and not upon the 
full joint law) is tailor made to the applications we have mind: As explained in the previous paragraph, we want to handle games
in which  
players interact with one another through the law of the control. In order to guess the impact this may have on the FBSDE representation, it is 
worth recalling that, in the problem 
    \eqref{SSCP},
 the optimal control may be represented in a quite generic way through the function $\hat{\alpha}$ with 
the adjoint variable $y$ therein taken as the gradient of the value function. 
Under the weak formulation approach, this turns the FBSDE system
\eqref{weak} into an FBSDE system depending on the marginal laws of $Z$, as $Z_{t}$ is known to have the representation 
$Z_{t} = \nabla_{x} u(t,X_{t}) \sigma$. Under the Pontryagin approach, it turns the FBSDE system 
\eqref{Pontryagin} into an FBSDE system depending upon the marginal laws of $Y$. Hence our choice to handle
systems of the form
\eqref{general_equation}. Still the reader should observe that, in order to handle the more general case when the players interact through the joint law
of the state and of the control, it is necessary to address systems of the same type
as \eqref{general_equation} but with the convention that $[X_{t},Y_{t},Z_{t}]$ is understood as the joint law of the triple $(X_{t},Y_{t},Z_{t})$; as we just mentioned, we do not address this level of generality in the report. 

In \eqref{general_equation}, the diffusion process $X$ is coupled to the diffusion process $(Y,Z)$ through the functions $B$ and $F$, representing the drift of the forward process and the driver of the backward process, respectively. The functions $B$ and $F$ are assumed to be Lipschitz in each of their arguments on $([0,T], \mathbb{R}^{3}, (\mathcal{P}_2(\mathbb{R}))^{3})$ and $G$ Lipschitz on $(\mathbb{R}, \mathcal{P}_2(\mathbb{R}))$, namely, for $(x,y,z,x',y',z) \in \mathbb{R}^6$ and $(\mu,\nu,\lambda,\mu',\nu',\lambda') \in (\mathcal{P}_2(\mathbb{R}))^6$, we have:
\begin{equation*}
    \begin{split}
        |B(t',x',y',z',\mu',\nu',\lambda') - B(t,x,y,z,\mu,\nu,\lambda)|\leq\ & K_B \bigl(|t'-t|+|x'-x|+|y'-y|+|z'-z| \\
        &+ W_{2}(\mu',\mu)+W_{2}(\nu',\nu) + W_{2}(\lambda',\lambda)\bigr) \\
        |F(t',x',y',z',\mu',\nu',\lambda') - F(t,x,y,z,\mu,\nu,\lambda)|\leq\ & K_F\bigl(|t'-t|+|x'-x|+|y'-y|+|z'-z| \\ &+W_{2}(\mu',\mu)+W_{2}(\nu',\nu)+W_{2}|(\lambda',\lambda)\bigr) \\
        |G(x',\mu') - B(x,\mu)| \leq\ & K_G\bigl(|x'-x|+W_{2}(\mu',\mu)\bigr).
    \end{split}
\end{equation*}

The goal of this project is to study numerical methods for solving this general FBSDE system. At some point in the report, we will 
relax the Lipschitz condition of $F$ in the variables $z$  and $\lambda$ and address an FBSDE with a quadratic driver $F$, 
see our examples in Section \ref{Examples}.

\section{Algorithms}\label{Algorithms}
\hspace{5mm} We implement two algorithms for numerically solving the FBSDE system in Equation (\ref{general_equation}). In the first algorithm, we represent paths of the stochastic processes $(X,Y,Z)$ using a tree structure, where branches of the tree represent quantization of the Brownian motion. In the second algorithm, we no longer represent the paths of the process, but the marginal laws of the process. In this case, the law of the process is discretized on a fixed temporal and spatial grid.

In both cases, the representation serves as a basis for a Picard scheme for approaching the solution. For the first algorithm, 
 the Picard scheme is implemented in the form of a global Picard scheme on the whole process; for the second one, iterations are just performed on 
 the marginal laws of the process. Although Picard's method sounds very natural, 
this strategy suffers, whatever the algorithm, from a serious drawback as Picard iterations for forward backward systems are only expected to converge in small time. Basically, the value of $T$ for which the algorithm will converge depends on the 
coupling strength of the system; we make this fact clear for the global method by showing how $T$ depends on the Lipschitz coefficients $K_B$, $K_F$ and $K_G$. In any case, bifurcations can be observed when $T$ is increased, or equivalently, when the coupling strength between the two equations is increased. For our convenience (since it can be costly to increase $T$), we will fix $T$ and explore the convergence of the algorithms as we vary the coupling strength.
This is one first step in our report: Compare how the two algorithms suffer from the coupling strength between 
the forward and backward equations.

The second key feature of our report is to use, for both algorithms, a continuation in time, which allows us to extend the value of the coupling parameter for which the algorithms converge. 

In the following sections, we detail our Picard approaches
for the two algorithms and the continuation in time method.

\subsection{Global Picard Iteration on a Small Time Interval}
\hspace{5mm} The main difficulty in numerically solving the FBSDE system is the fact that, not only the forward component $X=(X_t)_{0 \leq t \leq T}$ and backward component $(Y,Z)=(Y_t, Z_t)_{0 \leq t \leq T}$ are coupled, but also they run in opposite directions. Thus, neither equation can be solved independently of the other, seemingly requiring us to manage both time directions simultaneously. 
Several strategies are conceivable to sort out this issue. A first one is to make use of the \textit{decoupling field} of the system in order to work with one time direction instead of two (roughly speaking, the decoupling field is the value function $u$ in the weak approach and its derivative
$\nabla_{x} u$ in the Pontryagin one). We will investigate this method for our second algorithm;
its implementation is indeed pretty subtle in the mean field setting and it leads to the aforementioned 
\textit{Picard method on the marginal laws}. 
For our first algorithm, however, we limit ourselves to a \textit{brute force} approach. To decouple the equations, we propose a
\textit{global} Picard iteration scheme, whose definition is as follows. For the initial and terminal data of the problem ($\xi$ and $G$), we want to define a mapping $\Phi_{\xi,G}$ that will take the $j-1$ Picard iterate and produce the $j$ Picard iterate:

\begin{equation*}
        \Phi_{\xi,G}: (X^{j-1},Y^{j-1},Z^{j-1},[X^{j-1},Y^{j-1},Z^{j-1}]) \mapsto (X^{j},Y^{j},Z^{j},[X^{j},Y^{j},Z^{j}]).
\end{equation*}
We define the decoupled Picard scheme $\Phi_{\xi,G}$ as the following:

\begin{enumerate}
    \item First, solve
        \begin{equation*}
        \begin{split}
        dX^{j}_t &= B(t,X^{j-1}_t,Y^{j-1}_t,Z^{j-1}_t,[X^{j-1}_t,Y^{j-1}_t,Z^{j-1}_t])dt + \sigma dW_t \\
        X^j_0&=\xi \in L^2(\Omega,\mathcal{F}_0,\mathbb{P};\mathbb{R}), \\
        \end{split}
        \end{equation*}
        for $X^j$ which gives us $[X^j]$.
    \item Next, solve
        \begin{equation*}
        \begin{split}
        dY^j_t &= -F(t,X^{j}_t,Y^{j-1}_t,Z^{j-1}_t,[X^{j}_t,Y^{j-1}_t,Z^{j-1}_t])dt + Z^{j}_t dW_t \\
        Y^j_T &= G(X^{j}_T,[X^{j}_T]), \\
        \end{split}
        \end{equation*}
        for $Y^j$ and $Z^j$ which gives us $[Y^j]$ and $[Z^j]$.
    \item Return $(X^{j},Y^{j},Z^{j},[X^{j},Y^{j},Z^{j}])$.
\end{enumerate}
After initializing $(X^{0},Y^{0},Z^{0},[X^{0},Y^{0},Z^{0}])$, we can define a sequence by:
\begin{equation*}
(X^{j},Y^{j},Z^{j},[X^{j},Y^{j},Z^{j}])=\Phi_{\xi,G}(X^{j-1},Y^{j-1},Z^{j-1},[X^{j-1},Y^{j-1},Z^{j-1}]).
\end{equation*}
If the sequence $(X^j,Y^j,Z^j,[X^j,Y^j,Z^j])_{j=1,\dots}$ converges to some $(X,Y,Z, [X,Y,Z])$, then:
\begin{equation*}
(X,Y,Z,[X,Y,Z])=\Phi_{\xi,G}(X,Y,Z,[X,Y,Z]),
\end{equation*}
and thus, $(X,Y,Z)$ is a solution to the original FBSDE system in Equation (\ref{general_equation}).

Picard schemes such as this one are only guaranteed to converge for a small time horizon, $T$, depending on the Lipschitz coefficients of the functions $B$, $F$ and $G$ (and in fact not only the convergence of the Picard sequence but also the solvability of the equation 
may fail on an arbitrary time interval). We illustrate this idea through the following example
(the reader may find the general case in  
\cite{Delarue02}): Let the driver function $F = 0$ and define the common Lipschitz coefficient $K = \max (K_B, K_G)$. The FBSDE system becomes:
\begin{equation*}
    \begin{split}
        dX_t &= B(Y_t)dt + \sigma dW_t, \ X_0 = \xi \in \mathbb{R} \\
        dY_t &= Z_t dW_t, \ Y_T = G(X_T, [X_T]).
    \end{split}
\end{equation*}

We write the equation above in the integral form and take the expectation conditional to the filtration $\mathcal{F}_t$ in the backward equation:
\begin{equation*}
\begin{split}
        &X_t = \xi + \int_0^t B(Y_s)ds + \sigma W_t \\
        &Y_t = G(X_T, [X_T]) - \int_{t}^{T} Z_s dW_s, \quad \textrm{i.e.}, \ Y_t = \mathbb{E}(G(X_T, [X_T]) | \mathcal{F}_t).
\end{split}
\end{equation*}

Instead of one single $X$, 
let us now consider two (initial) processes denoted by $\hat{X}$ and $\tilde{X}$. After one Picard iteration, we arrive at $\hat{X}'$ and $\tilde{X}'$. The Picard iteration is given by:
\begin{equation*}
\begin{split}
    \hat{Y}_t &= \mathbb{E}(G(\hat{X}_T, [\hat{X}_T]) | \mathcal{F}_t) \\
    \hat{X}'_t &= \xi + \int_0^t B(\hat{Y}_s)ds + \sigma W_t,
\end{split}
\end{equation*}
and similarly for $\tilde{X}'$.
From the forward component, we have the following estimate for the difference between $\hat{X}'$ and $\tilde{X}'$:

\begin{equation*}
    \begin{split}
    & |\hat{X}'_t - \tilde{X}'_t|^2 \leq \biggl|\int_0^t (B(\hat{Y}_s) - B(\tilde{Y}_s)) ds\biggr|^2 \leq t K^2 \int_0^t |\hat{Y}_s - \tilde{Y}_s|^2 ds,
    \\
    & \mathbb{E} \Bigl[ \sup_{0 \leq t \leq T} |\hat{X}'_t - \tilde{X}'_t|^2 \Bigr] \leq T^2 K^2  \mathbb{E} 
    \Bigl[
    \sup_{0 \leq t \leq T} |\hat{Y}_t - \tilde{Y}_t|^2
    \Bigr].
    \end{split}
\end{equation*}
Then we write for the backward component:
\begin{equation*}
    \begin{split}
         \mathbb{E}\Bigl[ \sup_{0 \leq t \leq T} |\hat{Y}_t - \tilde{Y}_t|^2
        \Bigr] &= \mathbb{E} \Bigl[ \sup_{0 \leq t \leq T} | \mathbb{E} (G(\hat{X}_T, [\hat{X}_T]) - G(\tilde{X}_T, [\tilde{X}_T]) | \mathcal{F}_t) |^2 \Bigr]
        \\
        &\leq  4  \mathbb{E}\bigl[ |G(\hat{X}_T, [\hat{X}_T]) - G(\tilde{X}_T, [\tilde{X}_T])| ^2 \bigr]
        \\
        &\leq 8 K^2 \mathbb{E} \Bigl[ 
        |\hat{X}_T - \tilde{X}_T|^2 + W_{2}\bigl( [\hat{X}_T] , [\tilde{X}_T]\bigr)^2 \Bigr].
    \end{split}
\end{equation*}
In the second to last inequality, we used Doob's martingale inequality for the martingale term $\mathbb{E} (G(\hat{X}_T, [\hat{X}_T]) - G(\tilde{X}_T, [\tilde{X}_T]) | \mathcal{F}_t)_{0 \leq t \leq T}$. Also we used the fact that $
        W_{2} ([X_t], [\tilde{X}_t] )^2 \leq \mathbb{E} [ |X_t - \tilde{X}_t|^2 ]$.
        
Combining the inequalities above for both forward and backward components, we obtain the following estimate:
\begin{equation*}
    \begin{split}
        \mathbb{E}\Bigl[ \sup_{0 \leq t \leq T} |\hat{X}'_t - \tilde{X}'_t|^2 \Bigr] \leq 16 T^2 K^4  \mathbb{E}
        \bigl[ |\hat{X}_T - \tilde{X}_T|^2 \bigr] \leq 16 T^2 K^4   \mathbb{E}\Bigl[ \sup_{0 \leq t \leq T} |\hat{X}_t - \tilde{X}_t|^2
  \Bigr]. 
   \end{split}
\end{equation*}
And then,
\begin{equation*}
    \begin{split}
         \sup_{0 \leq t \leq T} W_{2}\bigl( [\hat{X}'_t], [\tilde{X}'_t] \bigr)^2  \leq \sup_{0 \leq t \leq T} \mathbb{E} \bigl[ |\hat{X}'_t - \tilde{X}'_t|^2 \bigr] \leq \mathbb{E} \Bigl[ \sup_{0 \leq t \leq T} |\hat{X}'_t - \tilde{X}'_t|^2
         \Bigr].
    \end{split}
\end{equation*}

Finally, when $16 T^2 K^4 < 1$, i.e. $T \leq 4/K^2$, the mapping $\Phi_{\xi,G}$ realizes a contraction on the forward component and the Picard iteration defined above will converge to the fixed point, providing a solution to the original FBSDE. Thus, we have illustrated that Picard schemes are only expected to converge in small time or for small coupling (i.e. smaller Lipschitz coefficients).

Keeping in mind that we eventually want to describe numerical schemes, we define a solver \textit{picard}, that will implement a finite number, $J_p \in \mathbb{Z}^+$, of Picard iterations. Thus, we define:
\vspace{5pt}

\noindent \hspace{5pt} $\bullet$ \textit{picard}$(\xi,G)$:
\begin{enumerate}
\item Initialize $X_{t}=\xi$, $Y_{t}=0$, and $Z_{t}=0, \forall t\in[0,T]$.
    \item For $1\leq j \leq J_p$
    \begin{enumerate}
        \item[] $(X,Y,Z,[X,Y,Z]) =\Phi_{\xi,G}(X,Y,Z,[X,Y,Z])$
    \end{enumerate}
\item Return $(X,Y,Z,[X,Y,Z])$.
\end{enumerate}

\subsection{Continuation in Time of the Global Method for Arbitrary Interval/Coupling}

\hspace{5mm}
Of course, we would like the mapping $\Phi_{\xi,G}$ to realize a contraction to make sure that the Picard iteration converges. As we just 
explained, this is the case when the forward and backward processes have a small enough coupling strength or for small time horizon $T$. But in practice, this is not always the case: It may happen that the FBSDE system is uniquely solvable, but that we are unable to prove that the Picard sequence converges. In order to overcome this issue, we follow the approach introduced in Chassagneux, Crisan, Delarue \cite{Solver}. Basically, the point is to divide the time interval into smaller intervals, called \textit{levels}, and to apply a Picard solver recursively between the levels.

To define the levels, we fix a time mesh: $\{0=T_0 < T_1, \dots,T_k,\dots T_{N_\ell-1}< T_{N_\ell}=T\}$. We would like to use the \textit{picard} solver introduced in the previous section to apply the Picard iteration on a given level. Notice that for an arbitrary level $[T_k,T_{k+1}]$, we do not know the initial condition for $X_{T_k}$ or the terminal condition for $Y_{T_{k+1}}$. Thus, our current approximation of these values will be inputs to the \textit{picard} solver in place of $\xi$ and $G$. We would also like to modify
the solver \textit{picard} to also take the current estimate of $(X,[X])$ so that we don't have to start from scratch every time we wish to use \textit{picard}. Thus, we define for level $k$:
\vspace{5pt}

\noindent \hspace{5pt} $\bullet$ \textit{picard}$[k](X,Y_{T_{k+1}})$:
\begin{enumerate}
\item Initialize $Y_{t}=0$, and $Z_{t}=0, \forall t\in[T_k,T_{k+1})$.
    \item For $1\leq j \leq J_p$
    \begin{enumerate}
        \item[] $(X,Y,Z,[X,Y,Z]) =\Phi_{X_{T_k},Y_{T_{k+1}}}(X,Y,Z,[X,Y,Z]))$
    \end{enumerate}
\item Return $(X,Y,Z,[X,Y,Z])$,
\end{enumerate}
where $X=(X_{t})_{T_k \leq t\leq T_{k+1}}$, and similarly for $Y$ and $Z$ and their laws. Note that $\Phi_{X_{T_k},Y_{T_{k+1}}}$ is the same as $\Phi_{\xi,G}$ defined earlier except the time horizon is $[T_k,T_{k+1}]$ instead of $[0,T]$ and the initial and terminal conditions are given by $X_{T_k}$ and $Y_{T_{k+1}}$ instead of $\xi$ and $G$, respectively.
In particular, pay attention that we no longer consider the terminal condition in the form of a mapping (like $G$) 
but in the form a random variable (like $Y_{T_{k+1}}$). Implicitly, this requires to store, from one step to another, the full random variable $Y_{T_{k+1}}$; hence our choice below to use a tree.

Now that we have a solver \textit{picard} which will implement the Picard iteration for a given level, next, we want to define a global solver to apply a continuation in time. The global solver, called \textit{solver}, is recursively defined as follows for some $J_s \in \mathbb{Z}^+$. For a given level $k$, define:
\vspace{5pt}

\noindent \hspace{5pt} $\bullet$ \textit{solver}$[k](X_{T_k},[X_{T_k}]):$
\begin{enumerate}
    \item Initialize $X_t=X_{T_k}$, $Y_t=0$ , and $Z_t=0$, $\forall t \in [T_k,T_{k+1}]$.
    \item For $1\leq j \leq J_s$
    \begin{enumerate}
        \item $(Y_{T_{k+1}},[Y_{T_{k+1}}])=$\textit{solver}$[k+1](X_{T_{k+1}},[X_{T_{k+1}}])$
        \item $(X,Y,Z,[X,Y,Z])=$\textit{picard}$[k](X,Y_{T_{k+1}})$
    \end{enumerate}
    \item Return $(Y_{T_k},[Y_{T_k}])$.
\end{enumerate}
As before, it is important to notice that the entry $X_{T_{k}}$ in \textit{solver}
is a random variable; $[X_{T_{k}}]$ is its law. Having the two in our notations is a bit redundant, but 
we feel it is more transparent for the reader. The break condition of the recursion is given by the terminal condition: 
\begin{center}
\textit{solver}$[N_\ell](X_{T_{N_\ell}},[X_{T_{N_\ell}}])=(Y_{T_{N_\ell}},[Y_{T_{N_\ell}}])=(G(X_{T_{N_\ell}},[X_{T_{N_\ell}}]),\mathcal{L}(G(X_{T_{N_\ell}},[X_{T_{N_\ell}}])))$
\end{center}

The goal of the continuation in time is for a Picard iteration scheme to converge even for large coupling parameters or large time horizon. We will see in Section \ref{Examples} that the continuation in time successfully achieves this goal for our benchmark examples.
We refer to 
\cite{Solver} for its theoretical analysis. 

Thus far, we have described our Picard approach and a continuation in time method. Now we need to provide a scheme for discretizing the Picard iteration mapping $\Phi_{\xi,G}$. In this report, we implement a tree algorithm; 
we also give a variant of it, in the form of a grid algorithm. For both algorithms, we consider the uniform time mesh with time step $h = T/N_t > 0$ with $N_\ell<N_t \in \mathbb{Z^+}$ and $ t_i=ih, i=0,...,N_t$. For convenience, we will assume the coarse time mesh used to define the levels, $\{0=T_0 < T_1, \dots,T_k,\dots T_{N_\ell-1}< T_{N_\ell}=T\}$, is a subset of the fine time mesh $\{0=t_0 < t_1, \dots, t_{N_{t-1}} <t_{N_t}=T\}$.

\subsection{Tree Algorithm for the Global Method}
\hspace{5mm} The first implementation of the Picard iteration, $\Phi_{\xi,G}$, is the tree algorithm presented in Chassagneux, Crisan, Delarue \cite{Solver}. We now provide a brief presentation of this algorithm.

\subsubsection{Time Discretization}
\hspace{5mm} Our first step in developing a discretization for the Picard iteration $\Phi_{\xi,G}$ is to discretize the problem in the time domain. We use the decoupled scheme derived in Chassagneux, Crisan, Delarue \cite{Solver}, and repeated below for convenience. As noted above, we consider the uniform time mesh with time step $h = T/N_t > 0$ with $N_\ell<N_t \in \mathbb{Z^+}$ and $ t_i=ih, i=0,...,N_t$.

Step 1) in defining $\Phi_{\xi,G}$ requires solving the forward equation for $X^j$. We use the classical Euler scheme:
\begin{equation*}
    \begin{split}
        X_{t_{i+1}}^j &=  X_{t_{i}}^j + h \,  B(t_i,X^{j-1}_{t_i},Y^{j-1}_{t_i},Z^{j-1}_{t_i},[X^{j-1}_{t_i},Y^{j-1}_{t_i},Z^{j-1}_{t_i}]) + \sigma \Delta W_i , \ X^j_0 = \xi.
    \end{split}
\end{equation*}
Note that when we calculate $X_{t_{i+1}}^j$, the value of $X_{t_{i}}^j$ and its law is known and could be substituted for $X_{t_{i}}^{j-1}$ and its law in the drift function, $B$.

Step 2) in defining $\Phi_{\xi,G}$ requires solving the backward equation for $Y^j$ and $Z^j$. We derive the discrete-time scheme to approximate the backward component, see e.g. \cite{bouchard2004discrete}.The $Y^j$ component at time $t_i$ in the backward scheme is obtained by taking the expectation conditional to $\mathcal{F}_{t_i}$, denoted as $\mathbb{E}_{t_i}$, of the backward equation between $t_i$ and $t_{i+1}$. Let $\Delta W_i = W_{t_{i+1}} - W_{t_{i}}$ denote the forward Brownian increment between $t_{i}$ and $t_{i+1}$. The driver function $F$ is approximated by its value at time $t_i$ and we use the fact that $F$ at $t_i$ and $Z_{t_i}$ are $\mathcal{F}_{t_i}$-measurable, and $\mathbb{E}_{t_i}(\Delta W_i) =\mathbb{E}_{t_i}(W_{t_{i+1}}-W_{t_{i}}) = 0$:

\begin{equation*}
\begin{split}
Y^j_{t_{i}} &= Y^j_{t_{i+1}} + \int_{t_i}^{t_{i+1}} F(t,X^{j}_t,Y^{j-1}_t,Z^{j-1}_t,[X^{j}_t,Y^{j-1}_t,Z^{j-1}_t])dt - \int_{t_i}^{t_{i+1}}Z^{j}_t dW_t \\
& \approx Y^{j}_{t_{i+1}} + h \ F(t_i,X^{j}_{t_i},Y^{j-1}_{t_i},Z^{j-1}_{t_i},[X^{j}_{t_i},Y^{j-1}_{t_i},Z^{j-1}_{t_i}]) - Z^{j}_{t_i} \Delta W_{i} \\
Y^j_{t_i} &= \mathbb{E}_{t_i} (Y^j_{t_{i+1}}) + h \ F(t_i,X^{j}_{t_{i}},Y^{j-1}_{t_{i}},Z^{j-1}_{t_{i}},[X^{j}_{t_{i}},Y^{j-1}_{t_{i}},Z^{j-1}_{t_{i}}]) \\
Y^j_{T} &= G(X^{j}_T,[X^{j}_T]). \\
\end{split}
\end{equation*}

As for the $Z^j$ component, we multiply the approximation of $Y^j_{t_i}$ by the Brownian increment $\Delta W_{i}$, taking the conditional expectation $\mathbb{E}_{t_i}$, and using $\mathbb{E}((\Delta W_i)^2) = h$. By noticing that our scheme will never make use of $Z^j_T$, we can simply set the terminal condition for the $Z^j$ component to $0$.
\begin{equation*}
\begin{split}
Z^j_{t_i} (\Delta W_{i})^2 &\approx Y^j_{t_{i+1}} \Delta W_{i} + (-Y^j_{t_{i}} + h \ F(t_i,X^{j}_{t_{i}},Y^j_{t_{i}},Z^{j-1}_{t_{i}},[X^{j}_{t_{i}},Y^j_{t_{i}},Z^{j-1}_{t_{i}}])) \Delta W_{i} \\
Z_{t_i}^j &= h^{-1} \mathbb{E}_{t_i}(Y_{t_{i+1}}^j \Delta W_{i}), \ Z_T^j = 0.
\end{split}
\end{equation*}

Putting this together, the time-discretized decoupled forward-backward scheme for Picard iteration of our general FBSDE system (\ref{general_equation}) is the following:

\begin{equation}
\label{scheme0}
\left\{
    \begin{aligned}
    & X_{t_{i+1}}^j =  X_{t_{i}}^j + h \ B(t_i,X^j_{t_i},Y^{j-1}_{t_i},Z^{j-1}_{t_i},[X^j_{t_i},Y^{j-1}_{t_i},Z^{j-1}_{t_i})] + \sigma \Delta W_i \\
    & X^j_0 = \xi, \\
    & Y^j_{t_i} = \mathbb{E}_{t_i} (Y^j_{t_{i+1}}) + h \ F(t_i,X^{j}_{t_{i}},Y^{j-1}_{t_{i}},Z^{j-1}_{t_{i}},[X^{j}_{t_{i}},Y^{j-1}_{t_{i}},Z^{j-1}_{t_{i}}]) \\
    & Y^j_{T} = G(X^{j-1}_T,[X^{j-1}_T])\\
    & Z_{t_i}^j = h^{-1} \mathbb{E}_{t_i}(Y_{t_{i+1}}^j \Delta W_i) \\
    & Z^j_T =0.
    \end{aligned}
    \right.
\end{equation}

Other variants of this forward-backward scheme are possible. For example, in the forward scheme we could change $j$ to $j-1$ and $t_{i}$ to $t_{i+1}$ on the right hand side. It is easy to observe that the forward-backward system is decoupled because of lagged Picard indices $j-1$ and $j$. Thus, given $(X^{j-1},Y^{j-1},Z^{j-1},[X^{j-1},Y^{j-1},Z^{j-1}])$ where $X^{j-1}=(X^{j-1}_{t_i})_{0 \leq i \leq N_t}$, and similarly for $Y$ and $Z$, we can solve the backward scheme autonomously and then the forward scheme to obtain $(X^{j},Y^{j},Z^{j},[X^{j},Y^{j},Z^{j}])$. We have not fully provided a discrete scheme for $\Phi_{\xi,G}$ yet, however, because for a given $t_i$, we have not discretized $(X_{t_i},Y_{t_i},Z_{t_i})$. This is the goal of the next section.

\subsubsection{Spatial Discretization via Tree Structure}
\hspace{5mm} The forward-backward decoupled scheme (\ref{scheme0}) above looks quite simple and explicit. However, it still presents some difficulties for the numerical computation. Firstly, it is difficult to compute the conditional expectation in the backward scheme. Secondly, it is non trivial to compute the law of $X_{t_{i+1}}$ forward in time. Even if there was no drift, the computation would involve the convolution of the law of $X_{t_{i}}$ and a Gaussian law of the Brownian increment. Ultimately, we will need a spatial dicretization.

The approach of this algorithm is to approximate the Brownian increments using a simple binomial approximation: $\Delta W_i = \pm \sqrt{h}$ with probability $1/2$. This gives rise to a binomial tree for the forward scheme. Each node on the tree at depth $i$ represents a value of $X_{t_i}$, and has two children nodes representing the two possible values of $X_{t_{i+1}}$ (the ``up $\uparrow$" and the ``down $\downarrow$" value), conditioned on the value of $X_{t_i}$. The two values are computed as follows:
\begin{equation}
    X^j_{t_{i+1}}(\uparrow \downarrow) = X_{t_{i}}^j + h \ B(t_i,X^j_{t_i},Y^j_{t_i},Z^j_{t_i},[X^j_{t_i},Y^j_{t_i},Z^j_{t_i}]) \pm \sigma \sqrt{h}.
    \label{binEuler}
\end{equation}

Suppose that we use $M$ points $x_1,..., x_{M}$ for the approximation of the law $\xi$ of the forward process at the initial time, i.e. $[X_0] = \xi \approx \sum_{k=1}^M p^0_k \delta_{x_k}(\cdot)$. Then we have $M$ parallel binomial trees at each Picard iteration. For Picard iterate $j$ and time $t_i$, the number of nodes at depth $i$ is $M \times 2^i$ with values of $(X^j_{t_i},Y^j_{t_i},Z^j_{t_i})$ saved on the nodes of the tree at depth $i$. The marginal law of each process at time $t_i$ can be determined by looking at all the values on the nodes at depth $i$. The backward scheme can be easily computed on the binomial tree. At the last time step, $T=t_{N_t}$, we have $Y^j_{T} = G(X^j_{T},[X^j_{T}])$ for each of the $M \times 2^{N_t}$ nodes. The conditional expectation in the backward scheme at $t_i$ is simply the average of the ``up" and ``down" branches at $t_{i+1}$.

To initialize the $j=0$ Picard iterate as in the definition of the solver \textit{picard}, we want to set $X_{t_i}=\xi, \forall i \in \{0, \dots, N_t\}$. This amounts to taking each initial value $x_k$ and initializing its entire tree to this value, meaning that $X_{t_i}^0=x_k$ for all nodes at depth $i$ and for all $i\in \{1,\dots,N_t\}$ of the $k$-th binomial tree. We then begin the Picard iteration by applying the mapping $\Phi_{\xi,G}$ as detailed above. Using the binomial tree and approximation of Brownian increments, the forward-backward decoupled scheme becomes fully implementable.

\subsection{Picard Iteration on the Marginal 
Laws: 
a Grid Algorithm}
\hspace{5mm} The complexity of the tree algorithm is exponential with respect to the number of time steps $N_t$, since the size of the binomial tree is of order $2^{N_t}$. The exponential complexity becomes problematic and makes continuation in time much slower when we deal with large time horizons. In order to reduce the size of the tree, a natural idea is to make some ``recombination'' of the binomial tree. But since the drift function $B$ depends on the value of the process, the two branches ``up-down'' and ``down-up'' from the same node at time $t_i$ will not coincide at time $t_{i+2}$, in general. Instead of recombination, we may fix a spatial grid of a controllable size that the binomial tree can be projected onto, in order to avoid exponential complexity. This will be 
the first ingredient of our grid algorithm, which is the main novelty in our report.

Inspired by the paper of Delarue and Menozzi \cite{Grid}, where the authors used a spatial grid for the approximation of FBSDEs without mean-field interaction, we make an intensive use of the notion of \textit{decoupling field}, which is the second key ingredient of our strategy. Indeed, using the representation result in Proposition 2.2 in \cite{PDEforUV}, we know that there exist deterministic feedback functions $(u,v): [0,T] \cross \mathbb{R} \cross \mathcal{P}_2(\mathbb{R}) \rightarrow \mathbb{R}$ with $u$ a solution to a nonlinear PDE (on the space of measures) such that for a solution $(X,Y,Z)$ to the general FBSDE system in Equation (\ref{general_equation}):
\begin{equation*}
    Y_{t}=u({t},X_{t},[X_{t}]) \quad \text{ and } \quad 
    Z_{t}=v({t},X_{t},[X_{t}]).
\end{equation*}
Generally speaking, $u$ is called the decoupling field of (\ref{general_equation}). Here comes the main observation: The time marginals of the solution to the general FBSDE system in Equation (\ref{general_equation}), $([X_t,Y_t,Z_t])_{0\leq t \leq T}$, can be equivalently characterized by $(\mu_t,u(t,\cdot),v(t,\cdot))_{0 \leq t \leq T}$ where  $[X_t]=\mu_t$. As in \cite{Grid}, our strategy is to thus approximate $(u,v)$ instead of $(X,Y,Z)$, but this is not so easy as 
$u$ and $v$ are defined on a space of infinite dimension (because the mean field component lives in 
${\mathcal P}({\mathbb R})$). In order to overcome this difficulty, we propose to freeze the mean field argument in 
$(u,v)$. This permits to regard $u$ and $v$ as functions of a finite-dimensional variable and thus to approximate both 
along the underlying spatial grid. Once an approximation of $u$ and $v$ has been computed for the given proxy of the marginal laws, 
we can update the value of this proxy by using a Picard method. Hence the name of this subsection.  
So, as opposed to the tree algorithm, we will no longer keep track of the pathwise laws of the processes. Instead, we will only compute the marginal laws at each time step of the time mesh by means of a Picard iteration. 

\subsubsection{Picard Iteration without Grid}

We first give the inputs and outputs of our new Picard mapping without any grid approximation:
\begin{equation*}
    \Psi_{[\xi],G}:(\mu^{j-1}_{t},u^{j-1}(t,\cdot),v^{j-1}(t,\cdot))_{0 \leq t \leq T} \mapsto (\mu^{j}_{t},u^{j}(t,\cdot),v^{j}(t,\cdot))_{0 \leq t \leq T},
\end{equation*}
with $\mu_0^{j-1}=\mu_0^{j}=[\xi]$. Denoting below $\varphi\sharp\nu$ as the push-forward measure of the measure $\nu$ by the function $\varphi$, $\Psi_{[\xi],G}((\mu^{j-1}_{t},u^{j-1}(t,\cdot),v^{j-1}(t,\cdot))_{0 \leq t \leq T})$ is defined by:

\begin{enumerate}
    \item Solve the following SDE for $X^j$:
        \begin{equation*}
        \begin{split}
        dX^{j}_t &= B(t,X^{j}_t,u^{j-1}(t,X^{j}_t),v^{j-1}(t,X^{j}_t),\mu^{j-1}_t,(u^{j-1}(t,\cdot),v^{j-1}(t,\cdot))\sharp \mu^{j-1}_t)dt + \sigma dW_t \\
        X^j_0&=\xi \in L^2(\Omega,\mathcal{F}_0,\mathbb{P};\mathbb{R}), \\
        \end{split}
        \end{equation*}
        Set $\mu^j_t:=[X^j_t]$.
    \item Next, find $(u^{j},v^j)(\cdot,\cdot) :[0,T] \cross \mathbb{R} \rightarrow \mathbb{R}^2$ through $(u^{j},v^j)(t,X^j_t):=(Y^j_t,Z^j_t)$ by solving:
        \begin{equation*}
        \begin{split}
        d Y^j_t &= -F(t,X^{j}_t,Y^{j}_t,Z^{j}_t,\mu^j_t,(u^{j-1}(t,\cdot),v^{j-1}(t,\cdot))\sharp \mu^{j}_t)dt 
        + Z^j_t dW_t \;, \\
        Y^j_T &= G(X^{j}_T,\mu^j_T)\;. \\
        \end{split}
        \end{equation*}
    \item Return $((\mu^{j}_{t},u^{j}(t,\cdot),v^{j}(t,\cdot))_{0 \leq t \leq T})$.
\end{enumerate}

Given $\Psi_{[\xi],G}$, we can go through the construction 
of the global Picard method and define similar routines \textit{picard} and \textit{solver}
for our new Picard approach
 by replacing formally 
$\Phi_{\xi,G}$ by $\Psi_{[\xi],G}$. The various inputs and outputs
in the new routines should be clear. 
The next step to make the whole algorithm entirely tractable is to show how to compute $\mu_{t}^j$, 
$u^j$, and $v^j$ explicitly in the mapping $\Psi_{[\xi],G}$. Similar to $\Phi_{\xi,G}$, this mapping will be defined explicitly for a discretized scheme on a temporal and spatial grid in the next two sections. As in the tree algorithm, we consider the uniform time mesh with time step $h = T/N_t > 0$ with $N_\ell<N_t \in \mathbb{Z^+}$ and $ t_i=ih, i=0,...,N_t$.

\subsubsection{Grid Approximation of Forward Component and its Law}
\hspace{5mm} We begin by fixing a spatial discretization grid. This grid could in principle be defined differently for each time step $t_i$, but for simplicity, we consider a homogeneous grid $\Gamma$ fixed for all time steps $t_i, \ i \in \{ 0, ..., N_t \}$ with constant spatial step size $\Delta x$. Let $\Pi$ be the projection function on the grid 
$$ \Gamma =  \{ x_k = x_1 + (k-1) \Delta x, \ k =  1, ..., N_x \}. $$ 
Precisely, $\Pi$ is given by
\begin{align*}
        \Pi(x) = x_k \text{ if } x \in [x-\Delta x/2,x+\Delta x/2) \quad \text{and} \quad \Pi(x) = x_0 \text{ if } x <x_0\,,\;
        					\Pi(x) = x_{N_x} \text{ if } x \ge x_{N_x}\,.
\end{align*}

The initial law $\xi$ of the forward process is approximated as $\xi \approx \mu_0(\cdot)$ on the grid $\Gamma$ with $N_x$ points. Recall that in the tree algorithm, we cannot choose $M$, the number of points for the approximation of the initial law, to be too large as we will need $M$ parallel binomial trees. Because the tree algorithm has exponential complexity, we are able to choose $N_x$ to be much larger than $M$, and in turn, the approximation of the initial law is more accurate in the grid algorithm than in the tree algorithm.
We can initialize the Picard iterate $(\mu^0_i,u^0_i,v^0_i)_{0 \leq i \leq N_t}$ similar as before by letting $\mu^0_{i}=\mu_0$ and $(u^0_i,v^0_i)=(0,0)$, for all $i \leq N_t$. 

We follow the definition of Step 1) for $\Psi_{[\xi],G}$,
but we use the Euler scheme for the forward process between $t_i$ and $t_{i+1}$:
\begin{equation} \label{eq euler scheme}
        X_{t_{i+1}}^{j} = X_{t_{i}}^{j} + h \ B(X_{t_i}^{j}, (u^{j-1}_i,v^{j-1}_i)(X_{t_{i}}^{j}), \mu^{j-1}_i,
        (u^{j-1}_i,v^{j-1}_i)\sharp\mu^{j-1}_i ) + \sigma \Delta W_{i}.
\end{equation}
Suppose the $j$ Picard iterate at time $t_i$ is given by $\mu^j_i$ with $\mu^j_0 = \mu_0$:
\begin{equation*}
\mu^j_i(\cdot) = \sum_{k = 1}^{N_x} p^j_{i,k} \delta_{x_k}(\cdot ), \ p^j_{i,k} \geq 0 \ \forall k \in \{ 1, ..., N_x \} \text{ and } \sum_{k = 1}^{N_x} p^j_{i,k} = 1.
\end{equation*}
The law of $X^j_{t_i}$ in the Euler scheme of the forward process is $\mu^j_i$. Then we would like to define $\mu^j_{i+1}$ as the law of $X_{t_{i+1}}^{j}$ in the Euler scheme above, but the quantity $X_{t_{i+1}}^{j}$ may not belong to the grid. Instead of using formula \eqref{eq euler scheme}, the natural idea is then to replace it by its projection on the grid $\Gamma$:
\begin{equation}
\label{euler_fwd} 
\begin{split}
    X_{t_{i+1}}^{j} &= \Pi \left(X_{t_{i}}^{j} + h \ B(X_{t_i}^{j}, (u^{j-1}_i,v^{j-1}_i)(X_{t_{i}}^{j}), \mu^{j-1}_i,
        (u^{j-1}_i,v^{j-1}_i)\sharp\mu^{j-1}_i ) + \sigma \Delta W_{i} \right) \\
    [X_{t_{i}}^{j}] &= \mu^j_i(\cdot) \text{ and } \mu^j_{i+1}(\cdot) = [X_{t_{i+1}}^{j}].
\end{split}
\end{equation}
The law of $X_{t_{i+1}}^{j}$ is the convolution of $\mu^j_{i}$ with transition density $q^j(t_i, t_{i+1}; x_k, x_n)$ with $x_k, x_n$ two points of the grid $\Gamma$ and $k,n \in \{ 1,..., N_x \}$. The transition probability is equal to the conditional probability that the process $X^j$ starting at time $t_i$ from point $x_k$ arrives at time $t_{i+1}$ at the point $x_n$, i.e. 

\begin{equation*}
q^j(t_i, t_{i+1}; x_k, x_n) = \mathbb{P} (X^{j}_{t_{i+1}} = x_n | X_{t_{i}}^{j} = x_k).
\end{equation*}
The discretized law $\mu^j_{i+1}$ is then written as:
\begin{equation*}
    \begin{split}
        \mu^j_{i+1}(\cdot) = (\mu_i * q^j(t_{i}, t_{i+1}))(\cdot) &= \sum_{n = 1}^{N_x} p^{j}_{i+1,n} \delta_{x_n}(\cdot) \\
        p^{j}_{i+1,n} = \sum_{k=1}^{N_x} \Bigl( p^{j}_{i,k} \times q^j(t_i, t_{i+1}; x_k, x_n) 
        \Bigr)
        &= \sum_{k=1}^{N_x} \Bigl( p^{j}_{i,k} \times \mathbb{P} (X^{j}_{t_{i+1}} = x_n | X_{t_{i}}^{j} = x_k) \Bigr).\\
    \end{split}
\end{equation*}
It is worth noticing that if we did not take the projection when computing $X_{t_{i+1}}^{j}$ in the scheme (\ref{euler_fwd}), then its law would be given by convolution of the law $\mu^j_{i}$ with the Gaussian transition density $\bar{q}^j(t_i, t_{i+1}; x_k, y)$ associated to the Euler scheme. The transition densities $q^j$ and $\bar{q}^j$ have the following relation:

\begin{equation*}
    \mathbb{P} (X^{j}_{t_{i+1}} = x_n | X_{t_{i}}^{j} = x_k) = q^j(t_i, t_{i+1}; x_k, x_n) = \int_{\beta(x_n, \Delta x/2)} \bar{q}^j(t_i, t_{i+1}; x_k, y) dy.
\end{equation*}

In fact, for a more tractable implementation of the forward scheme, we use the binomial approximation for the Brownian increments (\ref{binEuler}) introduced in the previous section. Note that quantization with more points can be easily applied. In this binomial case, with $(\uparrow)/(\downarrow)$ representing the ``up'' and ``down'' branches, respectively, the transition probabilities on the grid can be easily computed:
\begin{equation*}
    \begin{split}
        & \mathbb{P} (X^{j}_{t_{i+1}} = x_n | X_{t_{i}}^{j} = x_k) \\
        &=   \frac{1}{2} \left(\boldsymbol{1} (X^{j}_{t_{i+1}}(\uparrow) = x_n | X_{t_{i}}^{j} = x_k) + \boldsymbol{1}( X^{j}_{t_{i+1}}(\downarrow) = x_n | X_{t_{i}}^{j} = x_k )\right).
    \end{split}
\end{equation*}
Then we can write $\mu^j_{i+1}$ by computing the probabilities:

\begin{equation*}
\label{fwd_prob}
    p^j_{i+1,n} = \sum_{k=1}^{N_x} \frac{p^j_{i,k}}{2} \cdot \left(\boldsymbol{1} (X^{j}_{t_{i+1}}(\uparrow) = x_n | X_{t_{i}}^{j} = x_k) + \boldsymbol{1}( X^{j}_{t_{i+1}}(\downarrow) = x_n | X_{t_{i}}^{j} = x_k )\right)
\end{equation*}

At Picard iteration $j \geq 1$, the forward scheme finally gives the flow of measures $(\mu^j_i)_{i=0}^{N_t}$ at discrete time steps $(t_i)_{i=0}^{N_t}$ of the discretized law defined on the grid $\Gamma$. Thus, we have described an implementation of Step 1) in the definition of the Picard mapping $\Psi_{[\xi],G}$. In the next section, we detail the implementation of the backward components in Step 2).

\subsubsection{Grid for the approximation of the backward component}

\hspace{5mm} Given the current Picard iterates $(\mu^{j}_{i}$, $u^{j-1}_i(\cdot), v^{j-1}_i(\cdot))$, $i = 0, \cdots, N_t$, we would like to detail Step 2) in the definition of $\Psi_{[\xi],G}$. Since we have a discrete spatial grid $\Gamma$, we wish to compute the values of $u^j_i(x)$ and $v^j_i(x)$ for $x \in \Gamma$. 
By replacing $Y^j_{t_i}$ and $Z^j_{t_i}$ with their respective feedback functions in the backward component in Equation (\ref{scheme0}), for $x \in \Gamma $ we have the following backward scheme starting, for $i \le N_t -1$, with terminal condition at time $T = t_{N_t}$, $(u^j_{N_t},v^j_{N_t})=(G,0)$:
\begin{equation*}
    \begin{split}
        u^j_i(x) &= \mathbb{E} \left( u^j_{i+1}(X^{j}_{t_{i+1}},\mu^j_{i+1}) + h \cdot F(X^j_{t_{i}}, u^{j-1}_i(X_{t_{i}}^j), v^{j}_{i}(X_{t_{i}}^j),\mu^j_i, (u^{j-1}_i,v^{j-1}_i)\sharp \mu^j_i) 
        \, \vert \, X^j_{t_{i}}=x
        \right), \\
        v^j_i(x) &= \mathbb{E} \left( u^j_{i+1}( X^{j}_{t_{i+1}}) \cdot \Delta W_i/h
        \, \vert \, X^j_{t_{i}} = x
         \right).\\
    \end{split}
\end{equation*}
The variable $X^{j}_{t_{i+1}}$ with law $[X^{j}_{t_{i+1}}] = \mu^j_{i+1}$ is given by the forward scheme presented in the previous section with starting point $X^j_{t_i}=x$. Notice that by construction of the forward scheme, $X^{j}_{t_{i+1}} \in \Gamma$ and $supp(\mu^j_{i+1}) = \Gamma$ so $u^j_{i+1}(\cdot)$ has been calculated and saved at time $t_{i+1}$. 

We have a more explicit scheme in the case of binomial approximation of the Brownian increments, which is used for the numerical results for our examples, with $X^{j}_{t_{i+1}}(\uparrow)$ and $X^{j}_{t_{i+1}}(\downarrow)$ as defined in the forward scheme, always conditional to $X^j_{t_i} = x$ and given $\mu^j_i$ at time $t_i$:

\begin{equation*}
    \begin{split}
        u^j_i(x) &= \frac{1}{2} \left( u^j_{i+1}( X^{j}_{t_{i+1}}(\uparrow)) + u^j_{i+1}( X^{j}_{t_{i+1}}(\downarrow)) \right)  
         + h \cdot F(x, u^{j-1}_i(x), v^j_i(x),\mu^j_i, (u^{j-1}_i,v^{j-1}_i)\sharp \mu^j_i), \\
        v^j_i(x)
        &= \frac{h^{-1/2}}{2} \left( u^j_{i+1}(X^{j}_{t_{i+1}}(\uparrow)) - u^j_{i+1}( X^{j}_{t_{i+1}}(\downarrow)) \right).
    \end{split}
\end{equation*}

We have now described a scheme for Step 2) in the definition of $\Psi_{[\xi],G}$. Putting this together with the previous section, we have described a fully implementable scheme for $\Psi_{[\xi],G}$. Note that we can use this Picard iteration mapping to define the analogue of the solvers \textit{picard} and \textit{solver}. Importantly, we can also apply the continuation in time method to the grid algorithm as well. In the next section, we apply these two methods, the tree and grid algorithms, to five example problems.

\section{Examples}\label{Examples}

\hspace{5mm} We have collected five example problems to test the algorithms presented in Section \ref{Algorithms}.

\subsection{Linear Example}

\hspace{5mm} The first example is a linear model which comes directly from Chassagneux, Crisan, Delarue \cite{Solver}, in which they implemented the tree algorithm. The system of interest is the following:

\begin{equation*}
\begin{split}
    dX_t&=-\rho \mathbb{E}\left(Y_t\right) dt+\sigma dW_t \\
    X_0&=x_0 \\
    dY_t&=-a Y_t dt+Z_t dW_t \\
    Y_T&=X_T. \\
\end{split}
\end{equation*}

For this problem, the solution is known explicitly:

\begin{equation*}
Y_0=\frac{x_0e^{aT}}{1+\frac{\rho}{a}(e^{aT}-1)}
\end{equation*}

\begin{figure}[!htb]
\centering
\begin{subfigure}{.4\textwidth}
\includegraphics[scale=0.45]{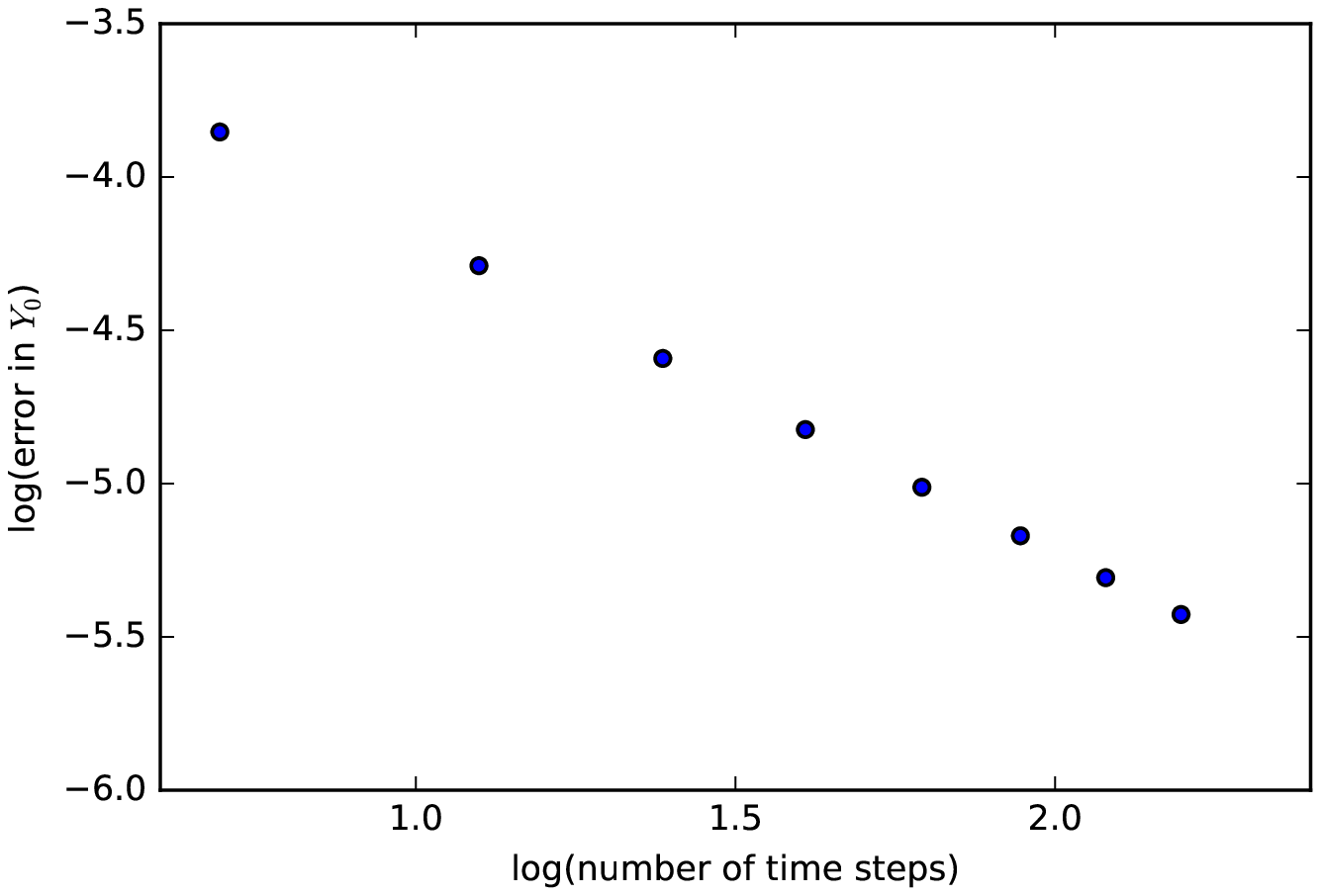}
\caption{Tree Algorithm}
\end{subfigure}
\begin{subfigure}{.4\textwidth}
\includegraphics[scale=0.45]{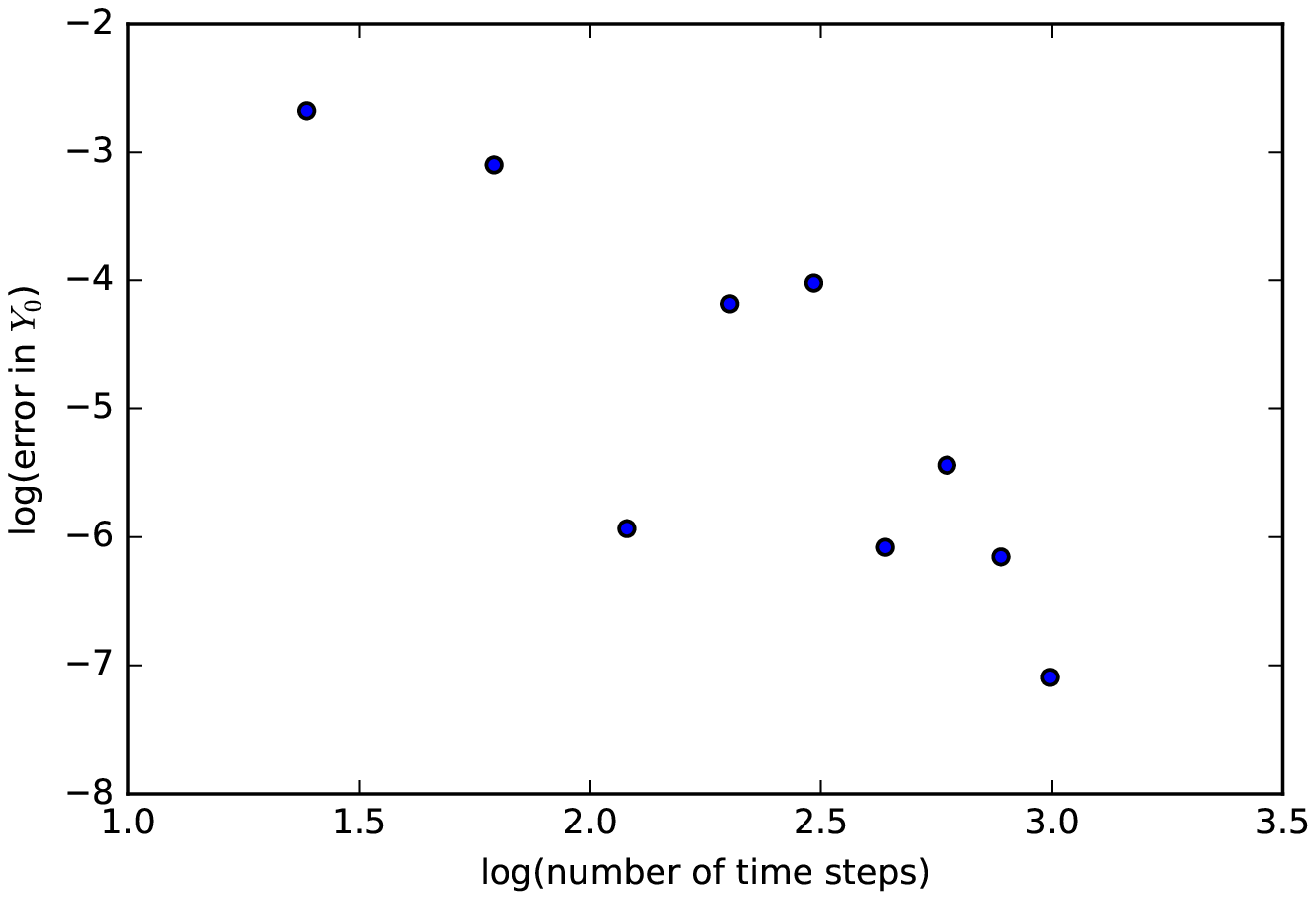}
\caption{Grid Algorithm}
\end{subfigure}
\caption{Linear Example: Convergence of the algorithms with one level as the number of time steps increases.}
\label{fig_1_1}
\end{figure}

For the numerical results, we let $\rho=0.1$, $a=0.25$, $\sigma=1$, $T=1$, and $x_0=2$. We vary $h$, the time step size, and for the grid algorithm, we use $\Delta x=h^2$. Figure \ref{fig_1_1} shows the log error between the numerical and true solution values of $Y_0$ as a function of log number of time steps for the tree and grid algorithm with one level (i.e. without using continuation in time). The left figure (tree algorithm) repeats the results in \cite{Solver} and as expected, it decreases linearly. On the other hand, we also observe a negative trend in the rate of convergence for the grid algorithm. Thus, both algorithms appear to converge to the true solution.

\subsection{Trigonometric Drivers Example}
\hspace{5mm} The second example also comes directly from Chassagneux, Crisan, Delarue \cite{Solver}. The system of interest is the following:

\begin{equation*}
\begin{split}
    dX_t&=\rho \cos\left(Y_t \right)dt+\sigma dW_t \\
    X_0&=x_0 \\
    Y_t&=\mathbb{E}_t\left(\sin(X_T)\right)
\end{split}
\end{equation*}

\begin{figure}[!htb]
\centering
\includegraphics[scale=0.5]{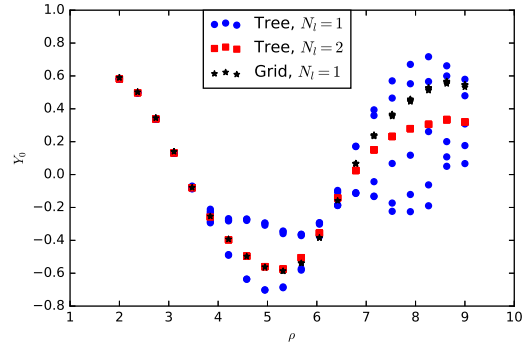}
\caption{Trigonometric Drivers Example: Bifurcations in the values of $Y_0$ appear as the coupling parameter, $\rho$, increases. The tree algorithm with one level is shown in blue circles and with two levels is shown in red squares. The grid algorithm with one level is shown in black asterisks.}
\label{fig_72_1}
\end{figure}

For the numerical results, we let $\sigma=1$, $T=1$, $x_0=0$, $h=1/6$ for the tree algorithm and $h=1/12$ and $\Delta x=h^2$ for the grid algorithm. For this problem, we observe a bifurcation when using the tree algorithm as we increase the coupling parameter, $\rho$. Figure \ref{fig_72_1} shows the values of $Y_0$ from the last $5$ Picard iterations. Starting at about $\rho=3.5$, the tree algorithm without continuation in time bifurcates. If the continuation in time method is used for the tree algorithm with two levels, there is no bifurcation for the range of values of $\rho$ showed in the plot. Note that the results from the tree algorithm repeat those in \cite{Solver}. The grid algorithm performs quite well in the sense that even with only one level, the algorithm converges for all of the values of $\rho$ in the plot. In particular, it avoids the exponential growth of the data structure characterizing the tree algorithm. Note that even though both the tree method with two levels and the grid method with one level converge, they produce different values for $Y_0$ for larger values of $\rho$. We believe this is because the tree algorithm is less accurate since the time step is larger.

\begin{figure}[!htb]
\centering
\includegraphics[scale=0.5]{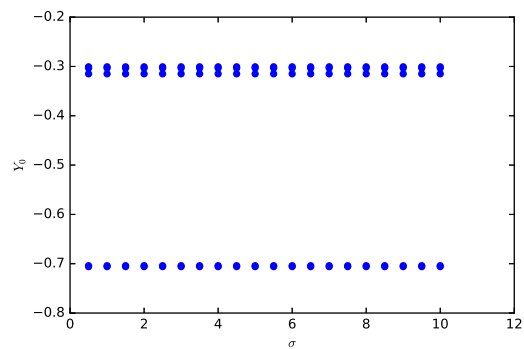}
\caption{Trigonometric Drivers Example: Tree algorithm with one level for $\rho=5$. Changing $\sigma$ has no effect on the bifurcation.}
\label{fig_72_2}
\end{figure}

\begin{figure}[!htb]
\centering
\begin{subfigure}{.4\textwidth}
\includegraphics[scale=0.45]{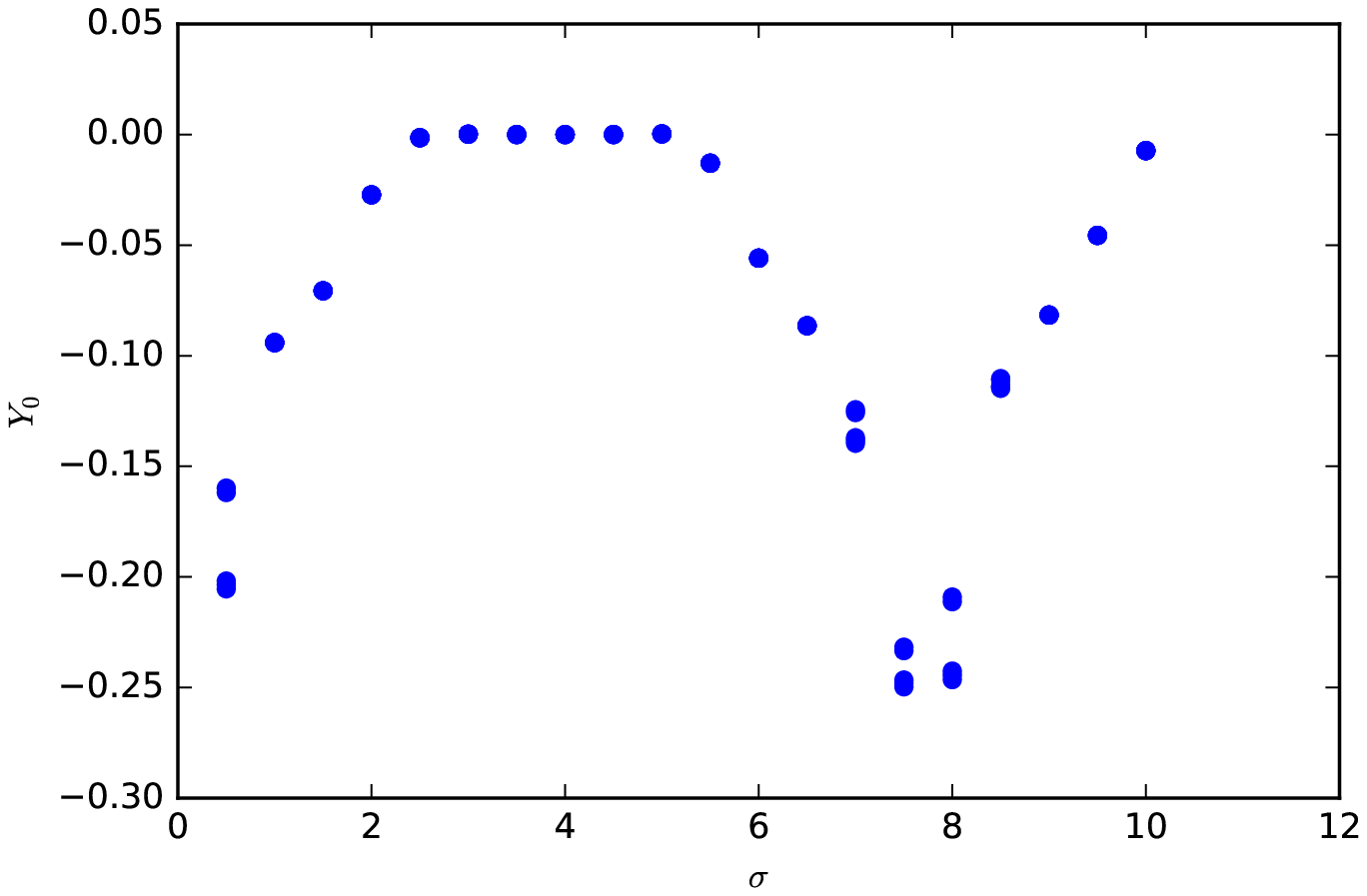}
\caption{$\rho=3.5$}
\end{subfigure}
\begin{subfigure}{.4\textwidth}
\includegraphics[scale=0.45]{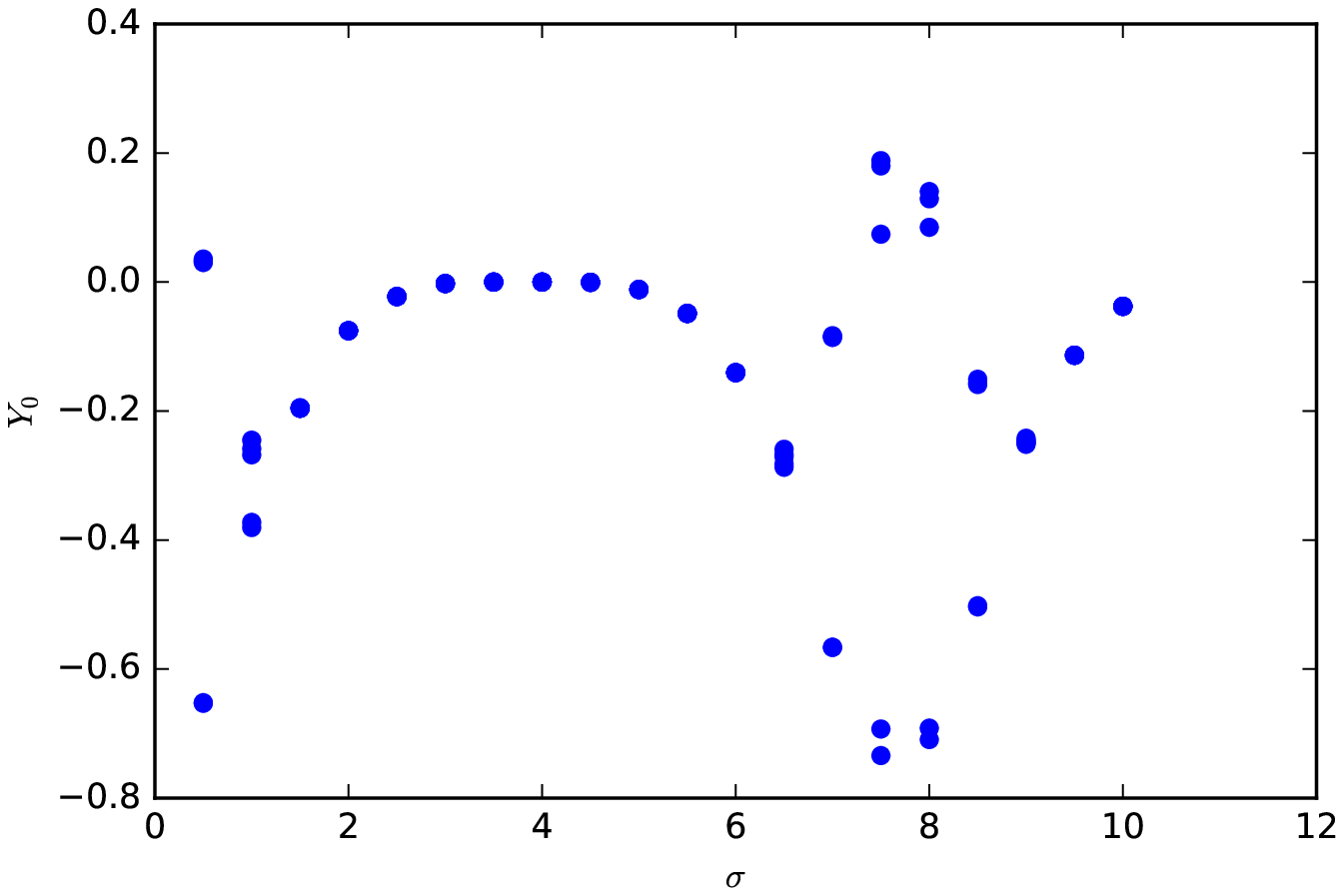}
\caption{$\rho=4$}
\end{subfigure}
\caption{Trigonometric Drivers Example: Tree algorithm with one level for $\rho=3.5$ and $\rho=4$. Changing $\sigma$ produces unexplained results.}
\label{fig_72_3}
\end{figure}

\begin{figure}[!htb]
\centering
\includegraphics[scale=0.5]{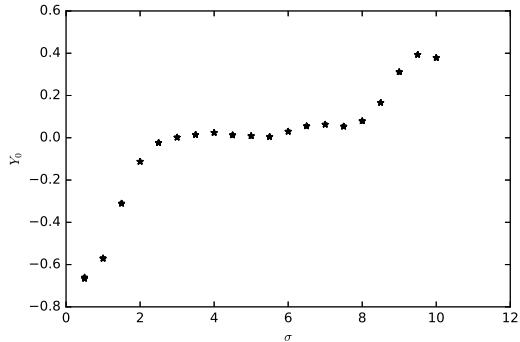}
\caption{Trigonometric Drivers Example: Grid algorithm with one level for $\rho=5$. There is no bifurcation as we vary $\sigma$.}
\label{fig_72_4}
\end{figure}

For values of $\rho$ for which the tree algorithm bifurcates, we were interested in the effect of changing $\sigma$ on the convergence. Figure \ref{fig_72_2} shows the last 5 values of $Y_0$ from the Picard iteration as $\sigma$ varies for $\rho=5$. Surprisingly, the value of $\sigma$ plays no role in the bifurcation for this value of $\rho$. To see if $\sigma$ would play a role when $\rho$ is closer to the bifurcation point, we have analogous plots when $\rho=3.5$ and $\rho=4$ (see Figure \ref{fig_72_3}). In these plots, however, it is clear that $\sigma$ does affect the bifurcation. Understanding the role of $\sigma$ on the bifurcation is an open question. It is interesting to note that even though the tree algorithm does not bifurcate for $\rho=3.5$ when $\sigma=1$, we still observe a bifurcation when we fix $\rho=3.5$ and vary $\sigma$. This suggests that we check if the grid algorithm bifurcates as we vary $\sigma$ for a fixed value of $\rho$ where there is no bifurcation when $\sigma=1$, such as $\rho=5$. Figure \ref{fig_72_4} shows that the grid algorithm does not bifurcate as we change $\sigma$ for fixed $\rho=5$.

\subsection{Mixed Model}\label{example_3}
\hspace{5mm} The third example also comes directly from Chassagneux, Crisan, Delarue \cite{Solver}. The system of interest is the following:

\begin{equation*}
\begin{split}
    dX_t&=- \rho Y_t dt+\sigma dW_t \\
    X_0&=x_0 \\
    dY_t&=\arctan \left( \mathbb{E}(X_t) \right) dt+Z_t dW_t \\
    Y_T&= \arctan \left( X_T \right).
\end{split}
\end{equation*}

\begin{figure}[!htb]
\centering
\includegraphics[scale=0.5]{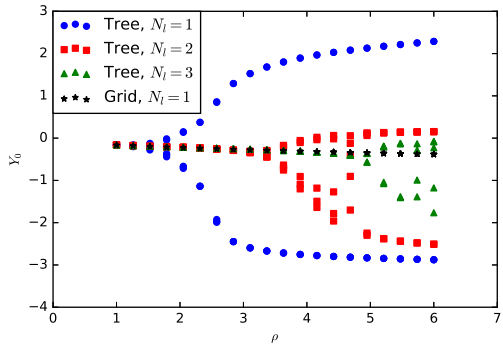}
\caption{Mixed Model: Bifurcations in the values of $Y_0$ appear as the coupling parameter, $\rho$, increases. The tree algorithm with one, two, and three levels is shown in blue circles, red squares, and green triangles, respectively. The grid algorithm with one level is shown in black asterisks.}
\label{fig_73_1}
\end{figure}

For the numerical results, we let $\sigma=1$, $T=1$, $x_0=2$, $h=1/6$ for the tree algorithm and $h=1/12$ and $\Delta x=h^2$ for the grid algorithm. We also observe a bifurcation for this problem as we increase the coupling parameter, $\rho$. Figure \ref{fig_73_1} shows the values of $Y_0$ from the last $5$ Picard iterations. Starting at about $\rho=1.5$, the tree algorithm without continuation in time bifurcates. If the continuation in time method is used for the tree algorithm with two levels, the bifurcation point is pushed back to about $\rho=3$, and pushed back further when using three levels to about $\rho=5$. Note that these results for the tree algorithm repeat those in \cite{Solver}. The grid algorithm performs quite well again in the sense that it converges for all of the values of $\rho$ shown, even when using only one level. Further, its major attractivness is the lower complexity compared to the tree algorithm. 

\begin{figure}[!htb]
\centering
\includegraphics[scale=0.5]{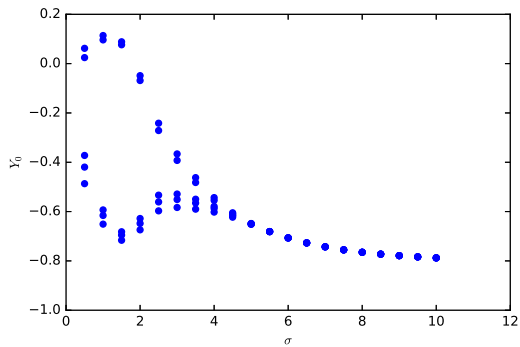}
\caption{Mixed Model: Tree algorithm with one level for $\rho=2$. The algorithm converges for large enough values of $\sigma$.}
\label{fig_73_2}
\end{figure}

As in the trigonometric drivers example, we can also investigate the effect of changing $\sigma$ for a value of $\rho$ where the tree algorithm without continuation in time bifurcates. For $\rho=2$, Figure \ref{fig_73_2} shows that the tree algorithm converges for large enough values of $\sigma$.

Since we have observed that the tree algorithm converges for small value of $\rho$ and large values of $\sigma$, this suggests trying a continuation in $\rho$ and/or $\sigma$, instead of the continuation in time. Instead of implementing a full continuation method, we used a simpler method of incrementing the parameter of interest.

\begin{figure}[!htb]
\centering
\includegraphics[scale=0.5]{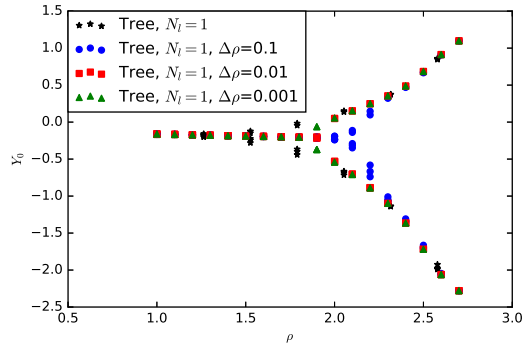}
\caption{Mixed Model: Tree algorithm with one level. Without continuation or incrementation is shown in black. Incrementation in $\rho$ with $\Delta \rho=0.1$, $0.01$, and $0.001$ are shown in blue circles, red squares, and green triangles, respectively.}
\label{fig_73_3}
\end{figure}

The incrementation in $\rho$ is performed by starting the algorithm with a small value of $\rho$ and letting it converge. Then $\rho$ is increased by some fixed $\Delta \rho$, and the algorithm is initialized with the solution from the previous value of $\rho$. Figure \ref{fig_73_3} shows the results from the incrementation in $\rho$. The incrementation in $\rho$ only increases the bifurcation point by a small amount.

\begin{figure}[!htb]
\centering
\includegraphics[scale=0.5]{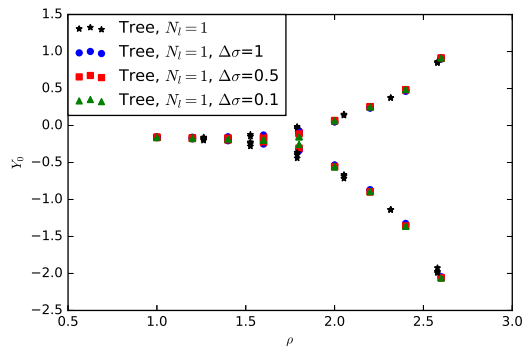}
\caption{Mixed Model: Tree algorithm with one level. Without continuation or incrementation is shown in black asterisks. Incrementation in $\sigma$ with $\Delta \sigma=1$, $0.5$, and $0.1$ are shown in blue circles, red squares, and green triangles, respectively.}
\label{fig_73_4}
\end{figure}

The incrementation in $\sigma$ is similar, except for each value of $\rho$, we start the algorithm with a sufficiently large value of $\sigma$ such that it will converge. Then $\sigma$ is decreased by a fixed $\Delta \sigma$. Figure \ref{fig_73_4} shows the results from the incrementation in $\sigma$. The incrementation in $\sigma$ also only increases the bifurcation point by a small amount.

If we take the previous example but change the drift and driver functions to also be in terms of the mean of the processes, then we have the following system, which we will refer to as the mixed model of means:

\begin{equation*}
\begin{split}
    dX_t&=- \rho \mathbb{E}(Y_t) dt+\sigma dW_t \\
    X_0&=x_0=2 \\
    dY_t&=\arctan \left( \mathbb{E}(X_t) \right) dt+Z_t dW_t \\
    Y_T&= \arctan \left(\mathbb{E} (X_T) \right).
\end{split}
\end{equation*}

\begin{figure}[!htb]
\centering
\includegraphics[scale=0.5]{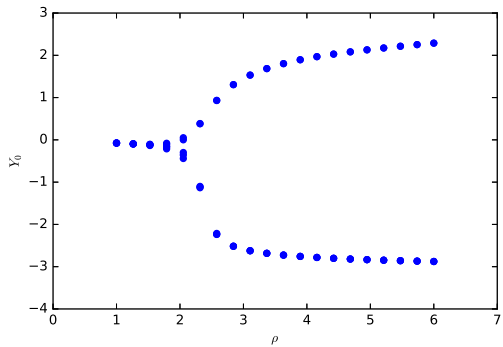}
\caption{Mixed Model of Means: Bifurcation for the tree algorithm with one level.}
\label{fig_73_E_1}
\end{figure}

\begin{figure}[!htb]
\centering
\includegraphics[scale=0.5]{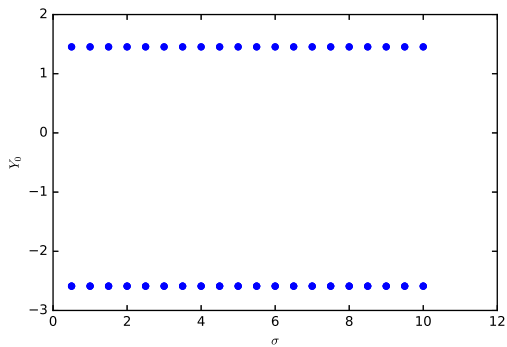}
\caption{Mixed Model of Means: Tree algorithm with one level with $\rho=3$. Changing $\sigma$ has no affect on the bifurcation.}
\label{fig_73_E_2}
\end{figure}

For the last example, we noticed that increasing $\sigma$ allows the tree algorithm to converge. We are thus interested to see if replacing the dynamics with the mean of the process will remove the effect of changing $\sigma$ on the convergence. We use the same values of the parameters as before. Figure \ref{fig_73_E_1} shows the bifurcation with one level of the tree algorithm (i.e. without using continuation in time). The effect of changing $\sigma$ is shown in Figure \ref{fig_73_E_2}. Our prediction is confirmed: $\sigma$ no longer affects the convergence when the dynamics are replaced with the mean of the process.

\subsection{Examples: Linear Quadratic Mean Field Games}\label{Truncation}

\hspace{5mm} The last two examples belong to the family of linear quadratic (LQ) games. In these models, the dynamics of the state are linear in the sense that the drift is defined by an affine function:
\begin{equation*}
b(t,x,\alpha)=A_t x + B_t \alpha + \beta_t.
\end{equation*}  
Furthermore, the running and the terminal cost are quadratic in the state and control variables. For the sake of simplicity, we choose not to include the cross terms, so that we define f and g as follows:
\begin{equation*}
\begin{split}
f(t,x,\alpha)&=\frac{1}{2} P_t x^2 + \frac{1}{2} Q_t \alpha^2 \\
g(x)&=\frac{1}{2} S x^2.  \\
\end{split}
\end{equation*}
The FBSDE system derived by approaching the LQ mean field games via the weak approach becomes:

\begin{equation}
    \begin{split}
    dX_t &= \big[A_t x + B_t \alpha + \beta_t\big] dt +
\sigma dW_t ,\quad  X_0=\xi \\
    dY_t &= \bigg[-\frac{1}{2} P_t X_t^2 + \big( A_t X_t + \beta_t\big) \frac{Z_t}{\sigma}+\frac{1}{2} \frac{B_t^2}{Q_t} \frac{Z_t^2}{\sigma^2}  \bigg] dt + dW_t,\quad Y_T = \frac{1}{2} S X_T^2.\\
    \end{split}
    \label{quadratic}
\end{equation}

We present two examples of linear quadratic games: a problem of flocking, and a price impact model of a trader.

\begin{remark}
We observe that the driver of the BSDE is quadratic in $Z_t$, which makes it seldom solvable. On the other hand, it is possible to obtain explicit solutions for LQ models. If we were to numerically solve equation (\ref{quadratic}), we might observe blow up from the quadratic terms in the backward equation and consequently the algorithm to not converge. Luckily, we do not observe such blowup in the examples we consider. But if blowup does occur, one could consider adapting the method presented by Chassagneux and Richou in \cite{Quadratic} to numerically approximate a quadratic BSDE.
\label{remark_trunc}
\end{remark}

\subsubsection{Flocking Problem}
\hspace{5mm} The next example problem models flocking. As in the paper of Nourian, Caines, and Malham\'{e} \cite{Nourian}, we consider the spatially homogeneous case where the state, $X_t$, represents the velocity of a representative player, or bird. Each bird controls its velocity through the process:

\begin{equation*}
dX_t=\alpha_t dt+\sigma dW_t,
\end{equation*}
where the control is chosen to minimize:
\begin{equation*}
J(\alpha)=\mathbb{E} \left[ \int_0^T \left(\frac{1}{2} \alpha_t^2+\frac{\rho}{2}(X_t-\bar{\mu}_t)^2 \right) dt\right].
\end{equation*}
over $\alpha \in \mathbb{A}$. Above, we let $\bar{\mu}_t$ denote the mean of the distribution $\mu_{t}$ (of the velocities of the birds) at time $t$. The running cost consists of two components. The first term encourages the birds to minimize their kinetic energy by not choosing a large control. The second term encourages the birds to align their velocities with the mean velocity of the group.

This model falls into the class of linear quadratic games. Assume that the initial condition is given by a constant, $X_0=x_0$. It can be shown that the solution is Gaussian with mean and variance:

\begin{equation*}
\begin{split}
\mathbb{E}(X_t)&=x_0 \\
Var(X_t)&=\sigma^2 \int_0^t \exp \left(-2 \int_s^t \eta_u du\right)ds,
\end{split}
\end{equation*}
where:
\begin{equation*}
\eta_t=\sqrt{\rho}\frac{e^{2\sqrt{\rho}(T-t)}-1}{e^{2\sqrt{\rho}(T-t)}+1}.
\end{equation*}
Using the weak formulation, the FBSDE system of interest is the following:

\begin{equation*}
\begin{split}
    dX_t&=-\frac{Z_t}{\sigma} dt+\sigma dW_t \\
    X_0&=x_0 \\
    dY_t&=-\left(\frac{Z_t^2}{2\sigma^2}+\frac{\rho}{2}(X_t-\mathbb{E}X_t)^2\right) dt+Z_t dW_t \\
    Y_T&=0.
\end{split}
\end{equation*}
If we use the Pontryagin formulation, the FBSDE system becomes:
\begin{equation*}
\begin{split}
    dX_t&=-Y_tdt+\sigma dW_t \\
    X_0&=x_0 \\
    dY_t&=-\rho\left(X_t-\mathbb{E}(X_t)\right) dt+Z_t dW_t \\
    Y_T&=0.
\end{split}
\end{equation*}

\begin{figure}[!htb]
\centering
\begin{subfigure}{.4\textwidth}
\includegraphics[scale=0.45]{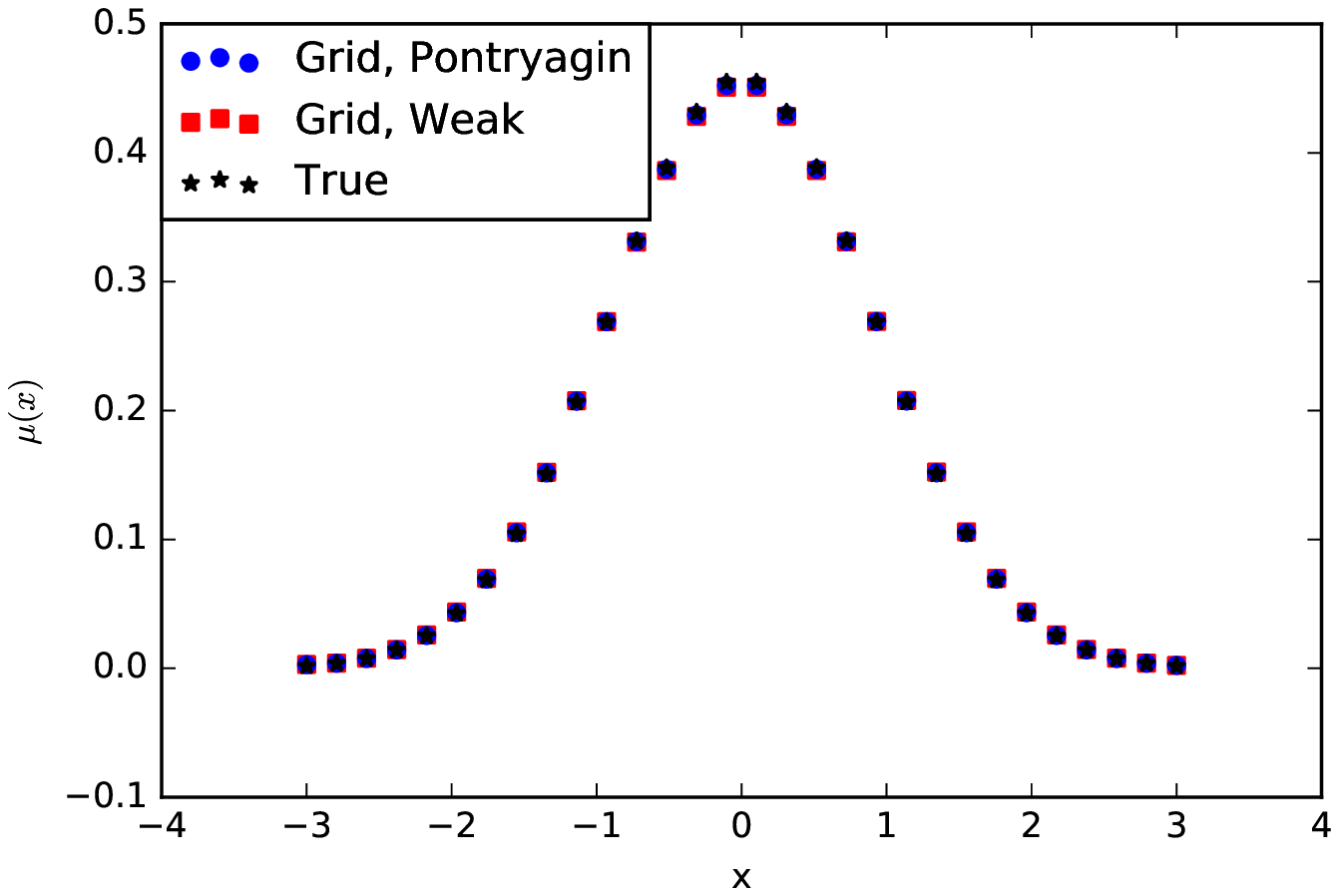}
\caption{}
\end{subfigure}
\begin{subfigure}{.4\textwidth}
\includegraphics[scale=0.45]{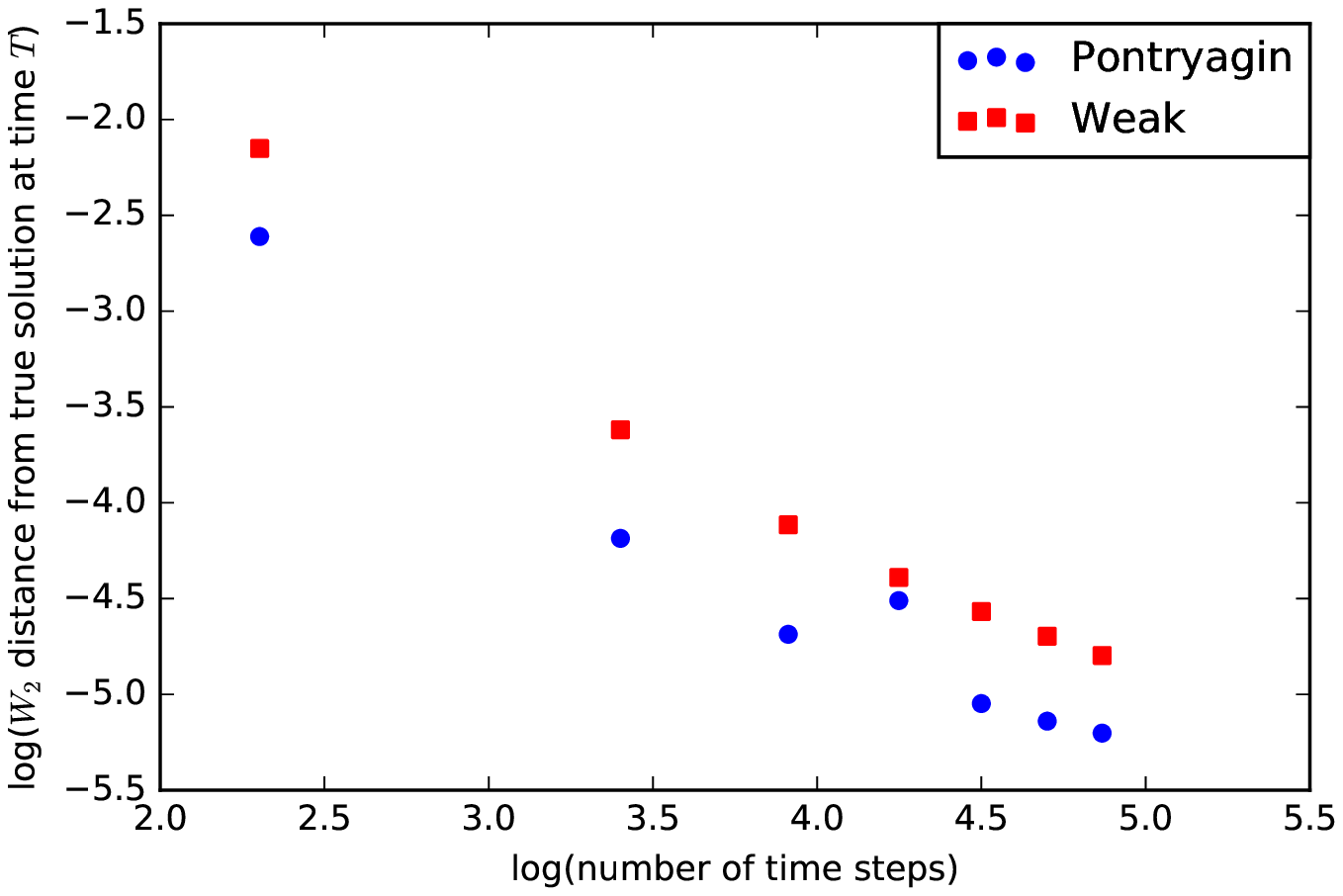}
\caption{}
\end{subfigure}
\caption{Flocking Problem: (a) Distribution $\mu_T$ of the players' states at time $T$ for the grid algorithm with one level. Pontryagin is in blue circles, weak is in red squares, true solution is shown in black asterisks.  (b) 2-Wasserstein distance between true solution and numerical solution for grid algorithm with one level as we increase the number of time steps, plotted as a log-log plot. Pontryagin approach is in blue circles and weak approach is in red squares.}
\label{fig_flocking_grid}
\end{figure}

The numerical results are presented for $\rho=1$, $\sigma=1$, $T=1$, $x_0=0$, $h=1/20$ for the tree algorithm and $h=1/130$ and $\Delta x=h^2$ for the grid algorithm. Figure \ref{fig_flocking_grid}(a) shows the results for the grid algorithm for both the weak and Pontryagin approaches. The plot shows the weights of the distribution $\mathcal{L}(X_{T})$. The results are similar between both approaches and coincide with the true solution.

We can look at the convergence rate by calculating the 2-Wasserstein distance between the numerical results and the true solution as we change the number of time steps. Since our state space is in one dimension, we can calculate the Wasserstein distance explicitly using the representation provided by Prokhorov \cite{W2}:

\begin{equation*}
    W_p(\mu,\nu)=\left(\int_0^1 | F_\mu^{-1}(u)-F_\nu^{-1}(u)|^pdu\right)^{1/p},
\end{equation*}
where $F_\mu(x)=\mu([0,x])$, denotes the cumulative distribution function. Figure \ref{fig_flocking_grid}(b) presents the convergence rate of the grid algorithm in terms of the 2-Wasserstein distance calculated between the true solution and numerical results with respect to the number of time steps. As expected, the 2-Wasserstein distance decreases towards 0 as we increase the number of time steps, for both the Pontryagin and weak approaches.

\begin{figure}[!htb]
\centering
\includegraphics[scale=0.5]{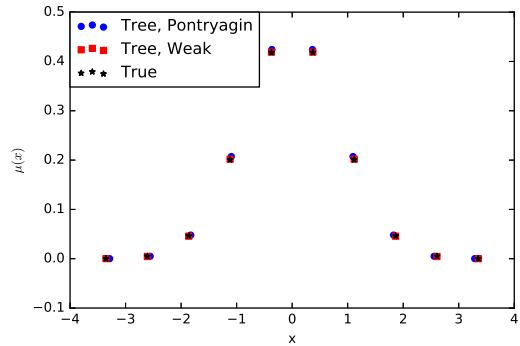}
\caption{Flocking Problem: Distribution $\mu_T$ of the players' states at time $T$ for the tree algorithm with one level. Pontryagin is shown in blue circles, weak is shown in red squares, and true solution is shown in black asterisks.}
\label{fig_flocking_tree}
\end{figure}

The results for the tree algorithm are shown in Figure \ref{fig_flocking_tree} for the weak and Pontryagin approaches. As with the grid algorithm, the weak and Pontryagin solutions are similar to each other and coincide with the true solution.

\subsubsection{Trader Problem}\label{trader}

\hspace{5mm} The last example shows an application of mean field games to finance, such as the trader congestion model in the paper of Cardialaguet and Lehalle, \cite{trader}. We focus on the \textit{Price Impact Model} presented in the book of Carmona and Delarue in \cite{MFGP}. The interest for this kind of model is motivated by their use in optimal execution problems for high frequency trading. Furthermore, it represents an instance of extended mean field game, also known as mean field game of controls, in which the representative agent interacts with the law of the control instead of the law of their state. The problem consists in a group of traders who have to buy or sell a large amount of shares in a given interval of time $[0,T]$. If they trade too fast, they will suffer from market impact. On the other hand, if they trade too slow, they will be affected by a large risk penalization. Approaching this problem as a mean field game, the inventory of the representative trader is modeled by a stochastic process $(X_t)_{0 \leq t \leq T} $ such that
\begin{equation*}
dX_t = \alpha_t dt +\sigma dW_t, \quad t\in [0,T],
\end{equation*}
where $\alpha_t$ corresponds to the trading rate. The price of the asset $(S_t)_{0 \leq t \leq T}$ is influenced by the trading strategies of all the traders trough the law of the controls $(\theta_t=\mathcal{L}(\alpha_t))_{0 \le t \le T}$ as follows:
\begin{equation*}
dS_t =\gamma \biggl( \int_{\mathbb{R}}a d\theta_t(a) \biggr) dt + \sigma_0 dW_t^0, \quad t\in [0,T],
\end{equation*} 
where $\gamma$ and $\sigma_0$ are constants and the Brownian motion $W^0$ is independent from $W$.
The amount of cash held by the trader at time $t$ is denoted by the process $(K_t)_{0 \le t \le T}$. The dynamic of $K$ is modeled by
\begin{equation*}
dK_t=-[\alpha_t S_t +c_{\alpha}(\alpha_t)]dt,
\end{equation*}
where the function $\alpha \mapsto c_{\alpha}(\alpha)$ is a non-negative convex function satisfying $c_{\alpha}(0)=0$, representing the cost for trading at rate $\alpha$. The wealth $V_t$ of the trader at time $t$ is defined as the sum of the cash held by the trader and the value of the inventory with respect to the price $S_t$:
\begin{equation*}
V_t=K_t+X_t S_t.
\end{equation*}
Applying the self-financing condition of Black-Scholes' theory, the changes over time of the wealth $V$ are given by the equation:
\begin{equation}
\begin{split}
    dV_t&=dK_t+X_tdS_t+S_t dX_t
    \\
   & =\Big[ -c_{\alpha}(\alpha_t)+\gamma X_t \int_{\mathbb{R}} a d\theta_t(a) \Big]dt + \sigma S_t dW_t + \sigma_0 X_t dW_t^0.
\end{split}
\label{wealth}
\end{equation}
We assume that the trader is subject to a running liquidation constraint modeled by a function $c_X$ of the shares they hold, and to a terminal liquidation constraint at maturity $T$ represented by a scalar function $g$. Thus, the cost function is defined by:
\begin{equation*}
J(\alpha)=\mathbb{E}\Big[\int_0^T c_X(X_t) dt +g(X_T) - V_T\Big].
\end{equation*}
Applying Equation (\ref{wealth}), it follows that
\begin{equation*}
J(\alpha)=\mathbb{E}\Big[ \int_0^T f(t,X_t,\theta_t,\alpha_t)dt +g(X_T)\Big],    
\end{equation*}
where the running cost is defined by
\begin{equation*}
f(t,x,\theta,\alpha)=c_{\alpha}(\alpha)+c_X(x)-\gamma x \int_{\mathbb{R}} a  d\theta(a),
\end{equation*}
for $0\leq t \leq T$, $x \in \mathbb{R}^d$, $\theta \in \mathcal{P}(A)$ and $\alpha \in A = \mathbb{R}$. We assume that the functions $c_X$ and $g$ are quadratic and that the function $c_{\alpha}$ is strongly convex in the sense that its second derivative is bounded away from $0$. Such a particular case is known as the Almgren-Chriss linear price impact model. Thus, the control is chosen to minimize:
\begin{equation*}
J(\alpha)=\mathbb{E}\left[ \int_0^T \left( \frac{c_{\alpha}}{2}{\alpha_t}^2+\frac{c_X}{2}X_t^2-\gamma X_t\int_{\mathbb{R}} a d\theta_t(a) \right)dt + \frac{c_g}{2}X_T^2\right],
\end{equation*}
over $\alpha \in \mathbb{A}$.
To summarize, the running cost consists of three components. The first term represents the cost for trading at rate $\alpha$. The second term takes into consideration the running liquidation constraint in order to penalize unwanted inventories. The third term defines the actual price impact. Finally, the terminal cost represents the terminal liquidation constraint.

As for the flocking example, this model falls in the class of linear quadratic games. Assume that the initial condition is given by a constant, $X_0=x_0$. The solution is Gaussian with mean and variance defined as:

\begin{equation*}
\begin{split}
\mathbb{E}(X_t)&=x_0 e^{-\int_0^t \frac{\bar{\eta}_s}{c_{\alpha}}ds}\\
Var(X_t)&=\sigma^2 \int_0^t e^{-\frac{2}{c_{\alpha}}\int_s^t \eta_r dr} ds\\
\end{split}
\end{equation*}
where:
 \begin{equation*}
\begin{split}
\bar{\eta}_t&= \frac{-C (e^{(\delta^+-\delta^-)(T-t)}-1)-c_g(\delta^+e^{(\delta^+-\delta^-)(T-t)}-\delta^-)}{(\delta^-e^{(\delta^+-\delta^-)(T-t)}-\delta^+)-c_gB(e^{(\delta^+-\delta^-)(T-t)}-1)} \\
\eta_t&=-c_{\alpha}\sqrt{c_X/c_\alpha}\frac{c_{\alpha}\sqrt{c_X/c_\alpha}-c_g-(c_{\alpha}\sqrt{c_X/c_\alpha}+c_g)e^{2\sqrt{c_X/c_\alpha}(T-t)}}{c_{\alpha}\sqrt{c_X/c_\alpha}-c_g+(c_{\alpha}\sqrt{c_X/c_\alpha}+c_g)e^{2\sqrt{c_X/c_\alpha}(T-t)}}, \\
\end{split}
\end{equation*}
for $t\in[0,T]$, where $B=1/c_{\alpha}$, $C=c_X, \delta^\pm=-D \pm \sqrt{R}$, 
with
$D = -\gamma /(2c_{\alpha})$
and
$R=D^2+BC$. \\
Using the weak approach yields the following FBSDEs system:
\begin{equation*}
    \begin{split}
    dX_t &= -\frac{1}{c_{\alpha}}\frac{Z_t}{\sigma} dt + \sigma dW_t,\quad X_0=x_0 \\
    dY_t &= -\left[\frac{c_X}{2} X_t^2 + \frac{\gamma }{c_{\alpha}}\frac{\mathbb{E}[Z_t]}{\sigma} X_t +\frac{1}{2c_{\alpha}} \left(\frac{Z_t}{\sigma}\right)^2\right] dt + Z_t dW_t,\quad Y_T = c_g \frac{X_T^2}{2}. \\
    \end{split}
\end{equation*}
Alternatively, the FBSDE system obtained via the Pontryagin approach is:
\begin{equation*}
\begin{split}
    dX_t &= -\frac{1}{c_{\alpha}}Y_t dt + \sigma dW_t, \quad X_0=x_0 \\
    dY_t &= -\left(c_X X_t + \frac{\gamma }{c_{\alpha}} \mathbb{E} [Y_t]\right) dt + Z_t dW_t,\quad Y_T = c_g X_T. \\
    \end{split}
\end{equation*}

\begin{figure}[!htb]
\centering
 \includegraphics[scale=0.5]{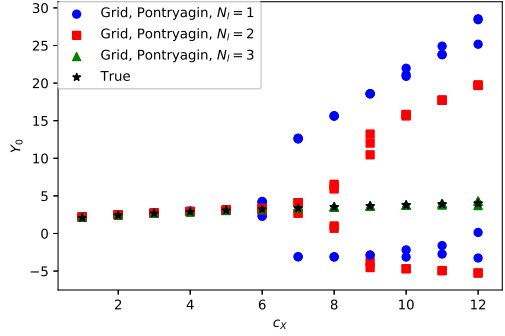}
\caption{Trader Problem: Bifurcations in the values of $Y_0$ depending on the coupling parameter $c_X$ for different number of levels in the grid algorithm. One, two, and three levels are shown in blue circles, red squares, and green triangles, respectively. The true value of $Y_0$ is shown in black asterisks.}
 \label{fig_trader_bif}
\end{figure}

The numerical results focus on the effect of the continuation method for the grid algorithm. In contrast with the previous examples, we show that the grid algorithm is also affected by bifurcation. Figure \ref{fig_trader_bif} shows the last five Picard iterations of $Y_0$ for the Pontryagin approach when the number of levels ranges from $1$ to $3$. Fixing the parameters $x_0=1$, $\sigma=0.7$, $1/c_{\alpha}=1.5$, $c_g=0.3$, $\gamma =2$, $T=1$, $h=1/12$, $\Delta x=h^2$ and increasing the coupling parameter, $c_X$, we observe that the bifurcation effect can be corrected by increasing the number of levels. In fact, Figure \ref{fig_trader_bif} shows that the true value of $Y_0$ matches the value computed numerically when using three levels.

\begin{figure}[!htb]
\centering
\begin{subfigure}{.4\textwidth}
\includegraphics[scale=0.45]{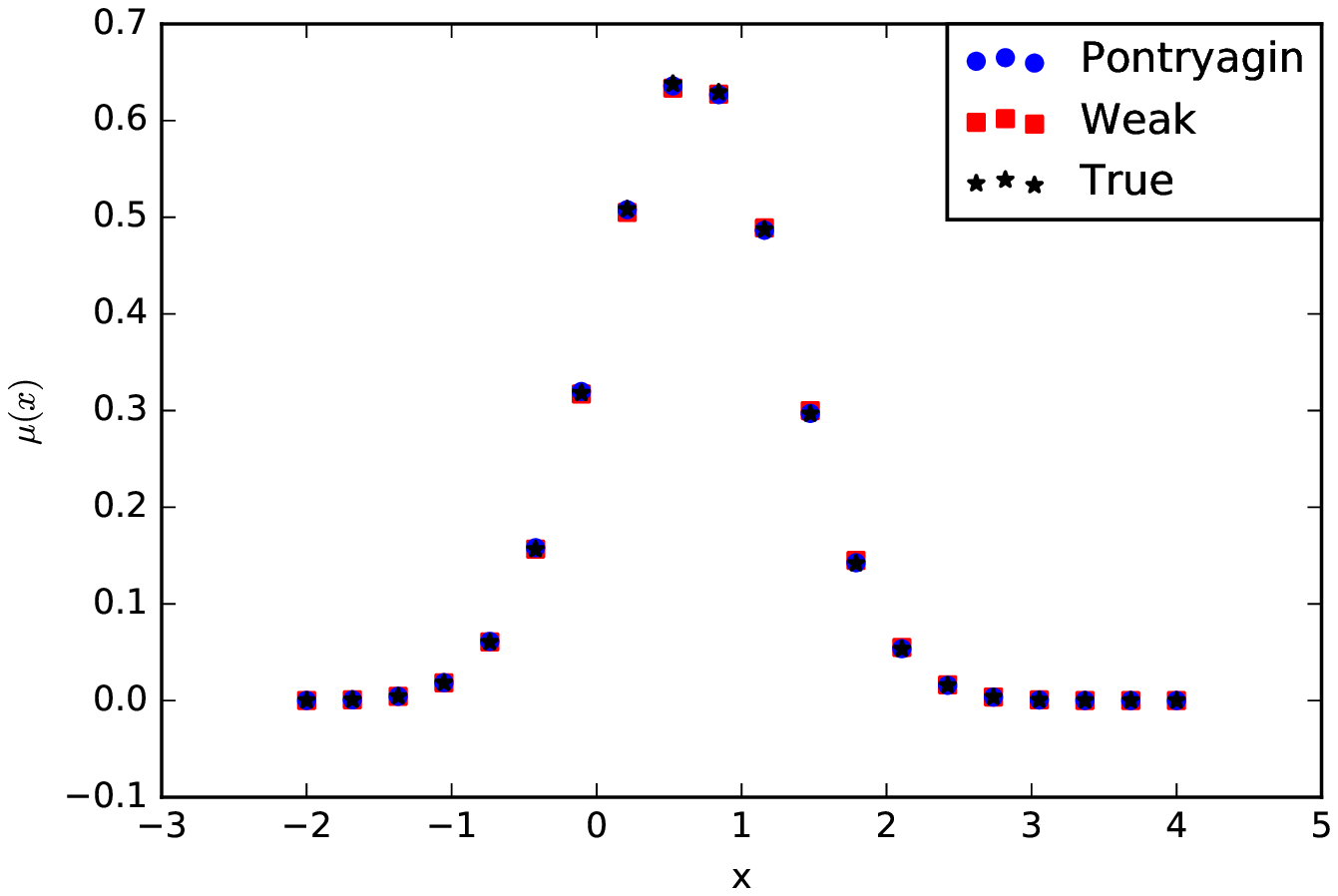}
\label{fig_trader_mu}
\caption{}
\end{subfigure}
\begin{subfigure}{.4\textwidth}
\includegraphics[scale=0.45]{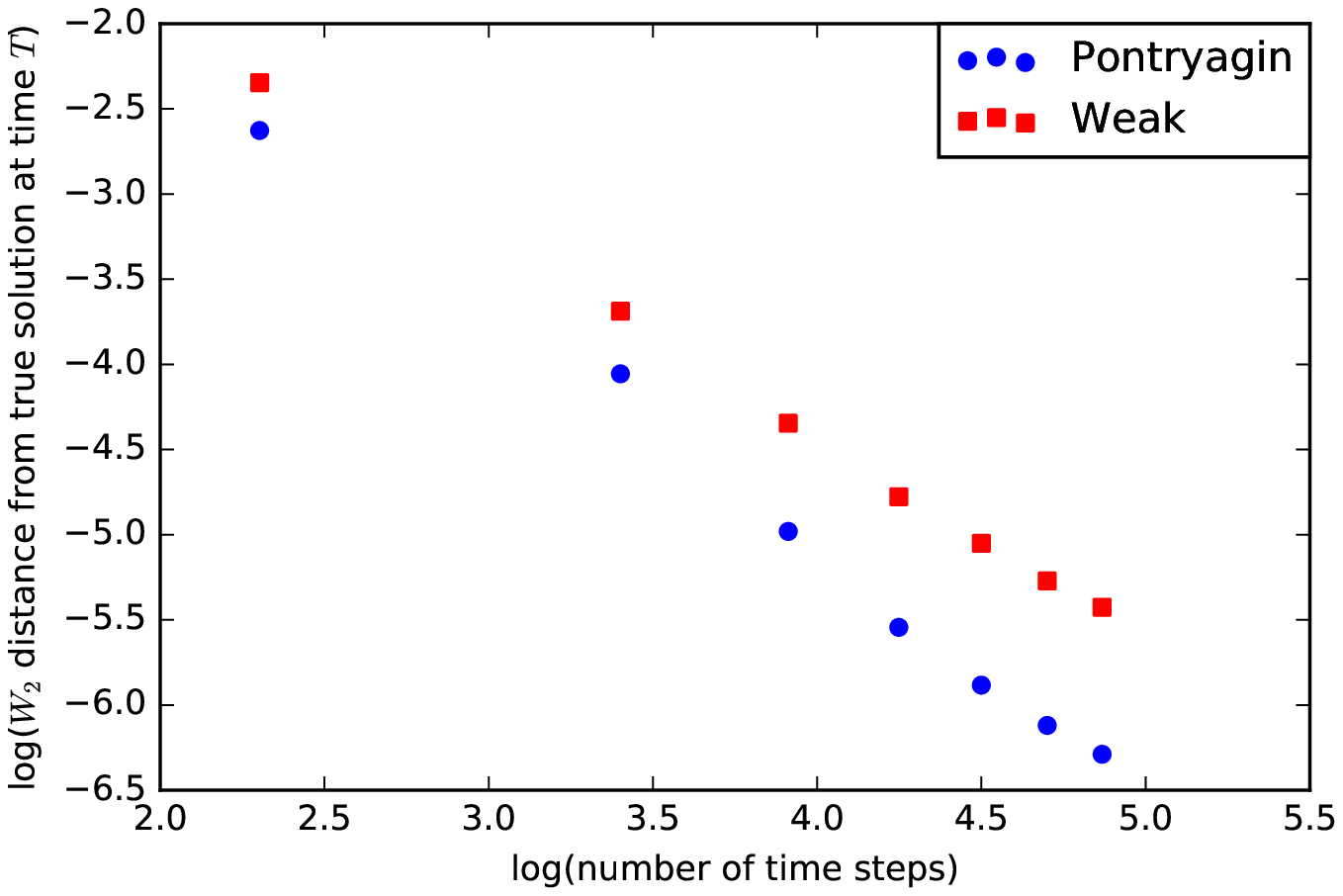}
\label{fig_trader_W2}
\caption{}
\end{subfigure}
\caption{Trader Problem: (a) Distribution $\mu_T$ of the players' states at time $T$ for the grid algorithm with one level. Pontryagin is shown in blue circles, weak is shown in red squares, and the true solution is shown in black asterisks. (b) 2-Wasserstein distance between true solution and numerical solution for grid algorithm with one level as we increase the number of time steps, plotted as a log-log plot. Pontryagin approach is shown in blue circles and weak approach is shown in red squares.}
\label{fig_trader_2}
\end{figure}

Furthermore, Figure \ref{fig_trader_2}(a) compares the distribution $\mathcal{L}(X_{T})$ obtained by the Pontryagin and the weak approaches using the grid algorithm with parameters $x_0=1$, $\sigma=0.7$, $1/c_{\alpha}=0.3$, $c_g=0.3$ $\gamma =2$, $c_X=2$, $T=1$, $h=1/130$ and $\Delta x=h^2$. The two approaches produce similar results that coincide with the true solution.

Figure \ref{fig_trader_2}(b) presents the convergence rate in terms of the 2-Wasserstein distance calculated between the true solution and numerical results with respect to the number of time steps. We again make use of the explicit representation of the Wasserstein distance \cite{W2}. The numerical solution is obtained using the grid algorithm with parameters $x_0=1$, $\sigma=0.7$, $1/c_{\alpha}=0.3$, $c_g=0.3$, $\gamma =2$, $c_X=2$, $T=1$, and $\Delta x=h^{2}$. As expected, the 2-Wasserstein distance decreases towards 0 as we increase the number of time steps.

\begin{figure}[!htb]
\centering
 \includegraphics[scale=0.5]{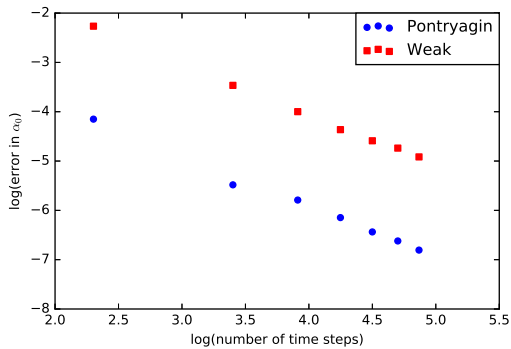}
\caption{Trader Problem: Error in the control at time 0, $\alpha_0$, as we increase the number of time steps, plotted as a log-log plot.}
 \label{fig_trader_error_control}
\end{figure}

The last plot, Figure \ref{fig_trader_error_control}, shows the error from the true solution of the control at time 0, $\alpha_0$, as we increase the number of time steps. This value is given by $\alpha_0=-Y_0/c_X$ for the Pontryagin approach and $\alpha_0=-Z_0/(c_X\sigma)$ for the weak approach. The true value is given by $\alpha_0=-\bar{\eta}_0 x_0/c_X$. As for the 2-Wasserstein distance, the error in the control decreases towards 0.

\section{Conclusion}\label{Conclusion}
\hspace{5mm} In conclusion, we have provided two algorithms for numerically solving FBSDEs of McKean-Vlasov type, which can be used to formulate the solutions to mean field game problems. The first algorithm is based on a pathwise tree structure. The second algorithm is based on a marginal grid structure. We have also proposed various refinements to the algorithms, including a continuation in time, and incrementation of a coupling parameter or the diffusion coefficient. The different numerical methods were illustrated on five benchmark examples.

The tree algorithm's main advantage is that we do not need to project the values of $X_t$ onto a discretized spatial grid, which potentially makes the algorithm more accurate. However, a significant disadvantage of the tree algorithm is the exponential growth of the data structure as the number of time steps is increased. This exponential growth is made worse yet if a higher order of quantization were to be used for approximating the Brownian increments.

The grid algorithm's main advantage is it avoids the exponential growth of the data structure. A higher order of quantization may be used without drastically changing the algorithm's complexity. A disadvantage of the grid algorithm is its sensitivity to the spatial step size with respect to the time step size. For the algorithm to be stable, the two step sizes need to be well adjusted to each other.

For both the tree and grid algorithms, we have observed that the continuation in time is able to extend the range of values of the coupling parameters for which the algorithms will converge. The incrementation methods proposed in Subsection \ref{example_3}, however, were not very successful at avoiding the bifurcations.

This report has touched on many things that could be explored deeper. First of all, the extension of the grid algorithm in \cite{Grid} to the mean field setting has not been studied from a theoretical standpoint. It is an open question to determine if this algorithm converges (meaning that the error decreases as the grid size decreases). The effect of the continuation in time or incrementation of the coupling parameter and/or diffusion coefficient has also yet to be studied. The numerical results also raised questions on the influence of the diffusion coefficient in the bifurcations.

\section{Acknowledgments}

\hspace{5mm} The authors would like to thank the organizers of the 2017 CEMRACS for the opportunity to study at CIRM. The student authors would also like to thank Professors Chassagneux and Delarue for mentoring us during the summer school. Finally, we would like to thank the funding sources that supported us during CEMRACS, including our respective universities, the NSF, and SMAI.

\end{document}